\documentclass[12pt]{article}
\usepackage{amsmath}
\usepackage{amssymb}
\usepackage{amsfonts}
\usepackage{eucal}
\usepackage{graphicx}
\numberwithin{equation}{section}
\newtheorem{Theorem}{Theorem}[section]
\newtheorem{Proposition}[Theorem]{Proposition}

\newtheorem{Remark}[Theorem]{Remark}
\begin{document}
\title{On the successive passage times of certain one-dimensional diffusions.}
\author{Mario Abundo$^*,$ Maria Beatrice Scioscia Santoro\thanks{Dipartimento di Matematica, Universit\`a  ``Tor Vergata'', via della Ricerca Scientifica, I-00133 Rome, Italy.
Corresponding author E-mail: \tt{abundo@mat.uniroma2.it}}
}
\date{}
\maketitle

\begin{abstract}
\noindent  We show in detail some results, outlined in a  previous paper regarding the case of Brownian motion (BM),
about the distribution of the $n$th-passage time of a one-dimensional diffusion obtained by a space or time transformation of
BM, through a constant barrier $a.$ Some explicit examples are reported.
\end{abstract}


\section{Introduction}
This is a continuation of \cite{abundo:phy16}, where we studied excursions of
Brownian motion with drift, and we found
explicitly the distribution of  the $n$th-passage time of Brownian motion (BM)
through a linear boundary.
In the present paper, by using the results of \cite{abundo:phy16}, we study in depth
the $n$th-passage time of a one-dimensional diffusion process $X(t),$ which is obtained
from BM by a space transformation or a time change. Particularly, focusing on
the second passage-time, we
find explicit formulae for the density of the second-passage time, and we calculate
its Laplace transform. We use
numerical computation and/or simulation, whenever analytical expressions cannot be achieved,
and we show graphically the results in a number of figures.
 \par
The successive-passage times of a diffusion $X(t)$ through a
boundary $S(t)$ are related to the excursions of $Y(t):=X(t) -
S(t);$ indeed, when $Y(t)$  is entirely positive or entirely
negative on the time interval $(s,u),$ it is said that it is an
excursion of $Y(t).$ Excursions have interesting applications in
Biology, Economics, and other applied sciences (see e.g.
\cite{abundo:phy16}, \cite{sch:39}). They are also related to the
last passage time of a diffusion through a boundary; actually,
last passage times of continuous martingales play an important
role in Finance, for instance, in models of  default risk (see
e.g. \cite{ell:2000}, \cite{jea:2009}). In the special case when
$X(t)$ is BM, its  excursions  follow the arcsine law, namely the
probability that BM has no zeros in the time interval $(s,u)$ is
$\frac 2 \pi \arcsin \sqrt {s/u} $ (see e.g. \cite{kleb:05}).
Really, by using Salminen's formula for the last passage time of
BM through a linear boundary (see \cite{salm:88}), in
\cite{abundo:phy16} we found the law of the excursions of drifted
BM, and  we derived the distribution of the $n$th passage time of
BM through a linear boundary; in this article, as we said, our aim is
to obtain analogous results for transformed BM through a constant boundary. \par The paper is organized as
follows: in Section 2, we recall the results of
\cite{abundo:phy16}, concerning the $n$th-passage time of BM
through a linear boundary, in Section 3 we extend these results to
certain diffusions $X(t),$ studying also the Laplace transform of
the second-passage time of $X(t).$ Section 4 reports some explicit examples,
while Section 5 is devote to conclusions and final remarks.

\section{Preliminary results on the nth-passage time of Brownian motion }
The first-passage time of BM  through the linear boundary $S(t)= a + bt,$ when starting
from $x ,$  is defined by $ \tau ^B_1 (S|x) = \inf \{ t >0: x +B_t = a + bt \} ,$ that will be denoted in this section  by
$\tau _1 (x),$ dropping, for simplicity,  the superscript $B$ which refers to BM and the dependence on $S;$  we recall
the Bachelier-Levy formula for the distribution of $\tau_1(x):$
$$
P(\tau_1 (x) \le t ) = 1 - \Phi((a -x) / \sqrt t + b \sqrt t )
 + exp(-2b (a-x)) \Phi(b \sqrt t -(a-x)/ \sqrt t ), \ t >0
$$
where $\Phi (y)= \int _ {- \infty } ^y \phi(z) dz ,$ with
$\phi (z) = e^ {- z^2/2}/ \sqrt {2 \pi },$
is the cumulative distribution function of the standard Gaussian variable.
If $(a -x)b > 0 ,$ then $P(\tau _1 (x) < \infty )= e^{-2b(a-x)},$ whereas, if
$(a -x)b \le 0 , \ \tau _1 (x)$  is finite with probability one and it has the Inverse Gaussian density,
which is non-defective
(see e.g. \cite{abundo:ric01}, \cite{kar:98}):
\begin{equation} \label{IGdensity}
f_ {\tau _1} (t) = f_ {\tau _1} (t |x) = \frac { d} {dt } P(\tau _1 (x) \le t) =  \frac {|a-x| } {t^ {3/2} } \ \phi \left ( \frac { a +bt -x } {\sqrt t } \right ) ,
\ t >0;
\end{equation}
moreover, if $b \neq 0,$ the expectation of $\tau _1 (x)$ is finite, being
$E(\tau _1(x)) = \frac {|a-x| } {|b| }.$ \par
The second-passage time of BM  through
$S(t),$ when starting from $x ,$ is defined by
$ \tau _2 (x) :=  \tau ^{B}_2 (S|x) = \inf \{ t > \tau _1(x): x +B_t = a + bt \} ,$ and generally, for $n \ge 2, \  \tau _n (x) = \inf
\{ t > \tau _{n-1}(x): x +B_t = a + bt \} $ denotes the $n$th-passage
time of BM through $S(t).$ \par
Now, we assume  that $ b \le 0 $ and $x < a ,$ or $b \ge 0 $ and $ x > a,$ so that  $P( \tau _1(x) < \infty )=1.$
For fixed $t>0,$ we consider the last-passage time prior to $t$ of BM, starting from $x ,$ through the boundary
$S(t)=a + bt,$
that is:
\begin{equation} \label{LPT}
\lambda _S ^t =
\begin{cases}
\sup \{ 0 \le u \le t : x + B_u = S(u) \}  &  {\rm if } \ \tau_1(x) \le t \\
0 & {\rm if } \ \tau_1(x) > t .
\end{cases}
\end{equation}
Notice that  \cite{abundo:phy16} contains a slight imprecision in the definition of $\lambda _S ^t ,$ and the function there denoted by $\psi _t (u),$ has
to be thought indeed as
the density  of $\lambda _S ^t ,$ conditional to the event $\{ \tau _1(x) \le t \};$ in other words, $\lambda _S ^t $ is zero with
probability $P( \tau _1 (x) >t),$ and it takes values spread on the interval $(0,t)$ with conditional density
$\psi _t (u)= P(\lambda _S ^t \in du | \tau _1(x) \le t ).$ However, this not affects the correctness of the results of \cite{abundo:phy16}.
In fact, it holds (see \cite{salm:88}, \cite{abundo:phy16}):
\begin{equation} \label{salminenEQ}
\psi _t (u)= \frac d {du} P (\lambda _S ^t \le u | \tau _1 (x) \le t ) =
\frac 1 { \sqrt {2 \pi u }} \exp (-b^2 u/2) \int _ {- \infty} ^{+ \infty }
\nu _ {x-a} (t-u, \widehat S )  dx, \ u \le t
\end{equation}
where $\widehat S (t) = S(t-u), \ u \le t ,$ and
\begin{equation} \label{nu}
\nu _x (v, \widehat S ) = \exp \{ - b(x-bt) - b^2 v /2 ) \} \frac {|x-bt|} {\sqrt { 2 \pi v^3} } \exp \{ - (x-bt)^2 /2v \}
\end{equation}
Then, the following explicit formula is obtained (see \cite{abundo:phy16}):
\begin{equation} \label{distrLPT}
 \psi _t (u)=  \frac { e^{ - \frac {b^2} 2 u} } { \pi \sqrt {u(t-u)}}
\left [ e^{ - \frac {b^2} 2 (t-u)} + \frac b 2 \sqrt { 2 \pi  (t-u)} \ \Big ( 2 \Phi ( b \sqrt {t-u} \ ) -1 \Big ) \right ]  , \ 0 <u <t
\end{equation}
Notice that
$ \psi _t $ is independent of $a;$
if $b=0,$ one gets
$
\psi _t (u) = \frac 1 { \pi \sqrt {u(t-u)}} \ ,  \ 0 <u <t ,
$
that is,  the arc-sine law with support in $(0,t).$
\begin{Remark}
For $z <t,$ the event $ \{ \lambda _S ^t \le z | \tau_1(x) \le t \} $ is nothing but the event $ \big \{ x+ B_u - S(u)$
{ \it has no zeros in the  interval} $(z,t) \big \} $.
\end{Remark}
In order to study the second-passage time, $\tau _2 (x),$ of $x+ B_t$ through the linear boundary $S(t)=a + bt ,$
when $ x <a $ and $b \le 0$ (the case when $b \ge 0$ and $x > a$ can be studied in a similar way),
we set $T_1(x)= \tau _1(x)$ and $T_2 (x) = \tau_2(x) - \tau _1 (x).$ It can be proved (see \cite{abundo:phy16}) that
$ T_2 (x)$ is finite with probability one, only if $b =0 .$
Conditionally to $\tau _1 (x)=s,$ the event $ \{ \tau _2 (x) > s+t \} \ (t >0) ,$ is nothing but the event $ \{ x + B_u - S(u)$
{ \it has no zeros in the interval } $(s, s+t) \} .$ Therefore, by using the above expression of $\psi _t (u)$ one obtains
$$ P( \tau _2 (x) > \tau _1 (x) +t | \tau _1 (x) =s ) =  P( \lambda _S ^{s+t} \le s | \tau _1 (x) \le s+t) = \int _0 ^s \psi _ {s+t} (y) dy $$
and so
\begin{equation} \label{conditionaldistrT2}
P( T_2(x) \le t | \tau _1 (x) =s ) = 1 - \int _ 0 ^s \psi _{s+t} (y) dy .
\end{equation}
Then:
\begin{equation} \label{distrT2}
P( T_2 (x) \le t ) = 1- \int _ 0 ^{+ \infty } f_{\tau _1} (s) ds \int _ 0 ^s \psi _{s+t} (y) dy .
\end{equation}
By taking the derivative with respect to $t,$ one gets the density of $T_2(x):$
\begin{equation} \label{denT2}
f_{T_2} (t) = \int _ 0 ^{+ \infty } f_{\tau _1} (s)  \left [ e^{ - b^2 (s+t)/2 } \frac {\sqrt s } {\pi \sqrt t (s+t) } \right ] \ ds .
\end{equation}
Notice that $ \int _ 0 ^s \psi _{s+t} (y) dy $ is decreasing in $t$
and
$ f_{T_2} (t) \sim  const / \sqrt t,$ as $ t \rightarrow 0^+.$
Moreover, from \eqref{distrT2} it follows that, if $b \neq  0,$ the distribution of $T_2(x)$ is  defective
In fact (see \cite{abundo:phy16}):
\begin{equation} \label{PT2infty}
P( T_2(x) = + \infty ) = 2 {\rm sgn} (b) E \left ( \Phi ( b \sqrt { \tau_1(x)} \ ) - \frac 1 2 \right ) >0,
\end{equation}
where ${\rm sgn} (b) =
\begin{cases}
0 & {\rm if} \  b=0 \\  |b| / b & {\rm if} \  b \neq 0
\end{cases}.$ \par\noindent
Since the function $g(s)= 2 {\rm sgn} (b) \left ( \Phi ( b \sqrt { s} \ ) - \frac 1 2 \right )$ is concave, by  Jensen's inequality
one gets
$ P( T_2(x) = + \infty ) \le \gamma (b),$ where $\gamma (b) = 2 {\rm sgn} (b) \left ( \Phi ( b \sqrt {  |(a-x)/b|  } \ ) - \frac 1 2 \right ) .$
Notice that  $ \gamma$ is an even function of $b.$ \par
On the contrary, if $b=0,$ from \eqref{PT2infty} we get that $P(T_2(x) = + \infty ) =0,$ that is, $T_2(x)$ is a proper random variable,
and
also $\tau _2 (x) = T_2(x) + \tau _1 (x)$ is finite with probability one;
precisely, by calculating the integral in \eqref{distrT2}
we have:
\begin{equation}
P( T_2(x) \le t ) = \int _0 ^{+ \infty } \frac 2 \pi \arccos \sqrt { \frac s {s+t} } \
\frac {|a-x| } {\sqrt { 2 \pi } s^{3/2} } \ e^{ - (a -x)^2 /2s } ds
\end{equation}
then, taking the derivative with respect to $t,$ the density of $T_2(x)$ turns out to be:
\begin{equation} \label{denT2b0}
f_{T_2} (t) = \int _ 0 ^{+ \infty } \frac {1} {\pi (s+t) \sqrt t }
\frac {|a-x| } {\sqrt { 2 \pi } s } \ e^{ - (a -x)^2 /2s } ds .
\end{equation}
In the Figure 1, taken from  \cite{abundo:phy16},
the plot of $P(\tau _2(x) = + \infty) = P(T _2(x) = + \infty),$
as a function of $b\le 0,$  is compared with the plot of its upper bound $\gamma (b),$  for $a=1$ and $x=0.$

\begin{figure}
\centering
\includegraphics[height=0.33 \textheight]{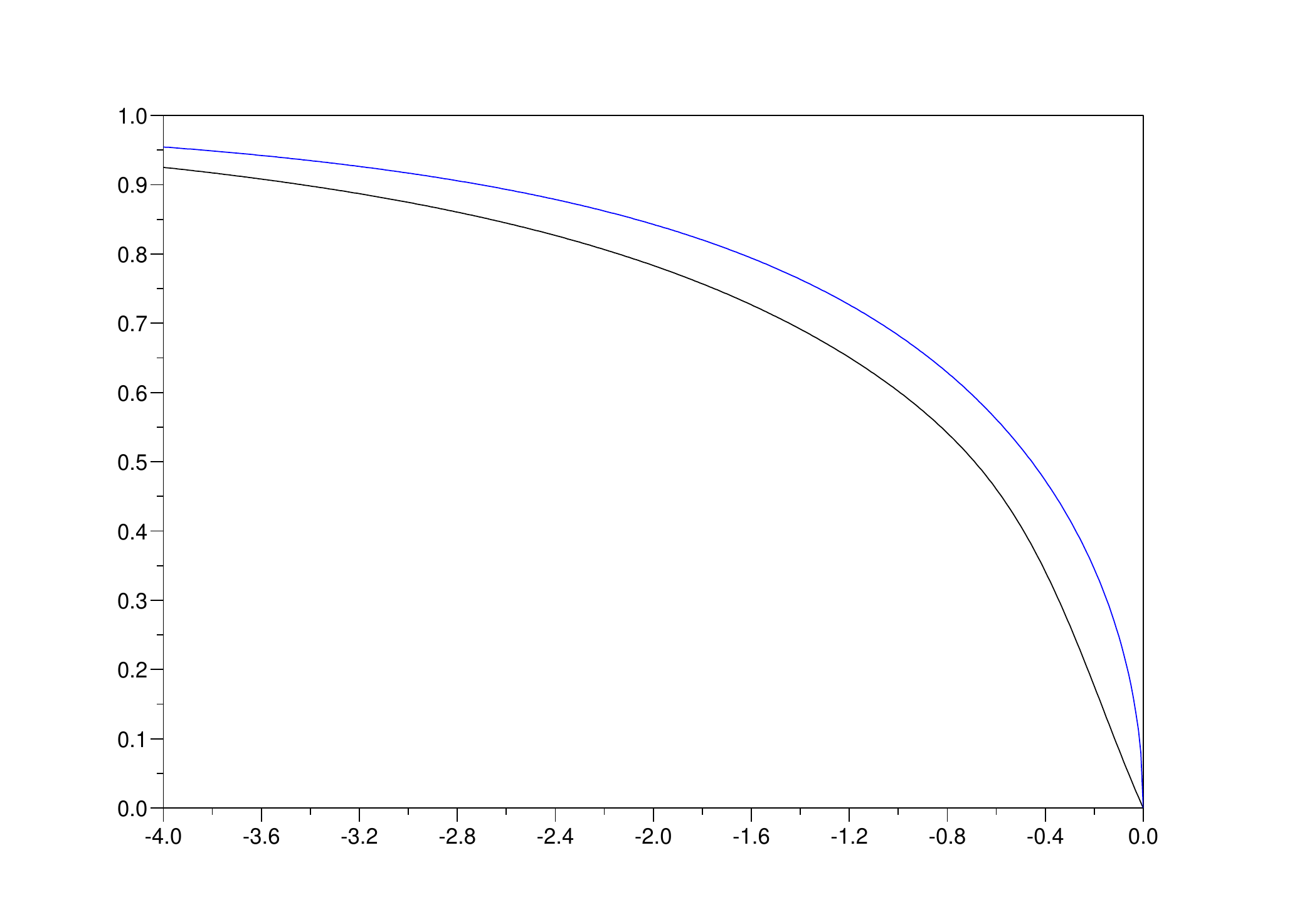}
\caption{Comparison of the shapes of $P(\tau _2(x) = + \infty)$ (lower curve)
and that of its upper bound $\gamma (b),$   as a function of $b\le 0,$ for $a=1$ and $x=0.$
}
\end{figure}

In the Figure 2, taken from \cite{abundo:phy16}, it is reported the probability density of $T_2(x),$ obtained from \eqref{denT2},
by calculating numerically the integral, for
$x=0, \ a=1,$ and various values of the parameter $b \le 0.$

\begin{figure}
\centering
\includegraphics[height=0.33 \textheight]{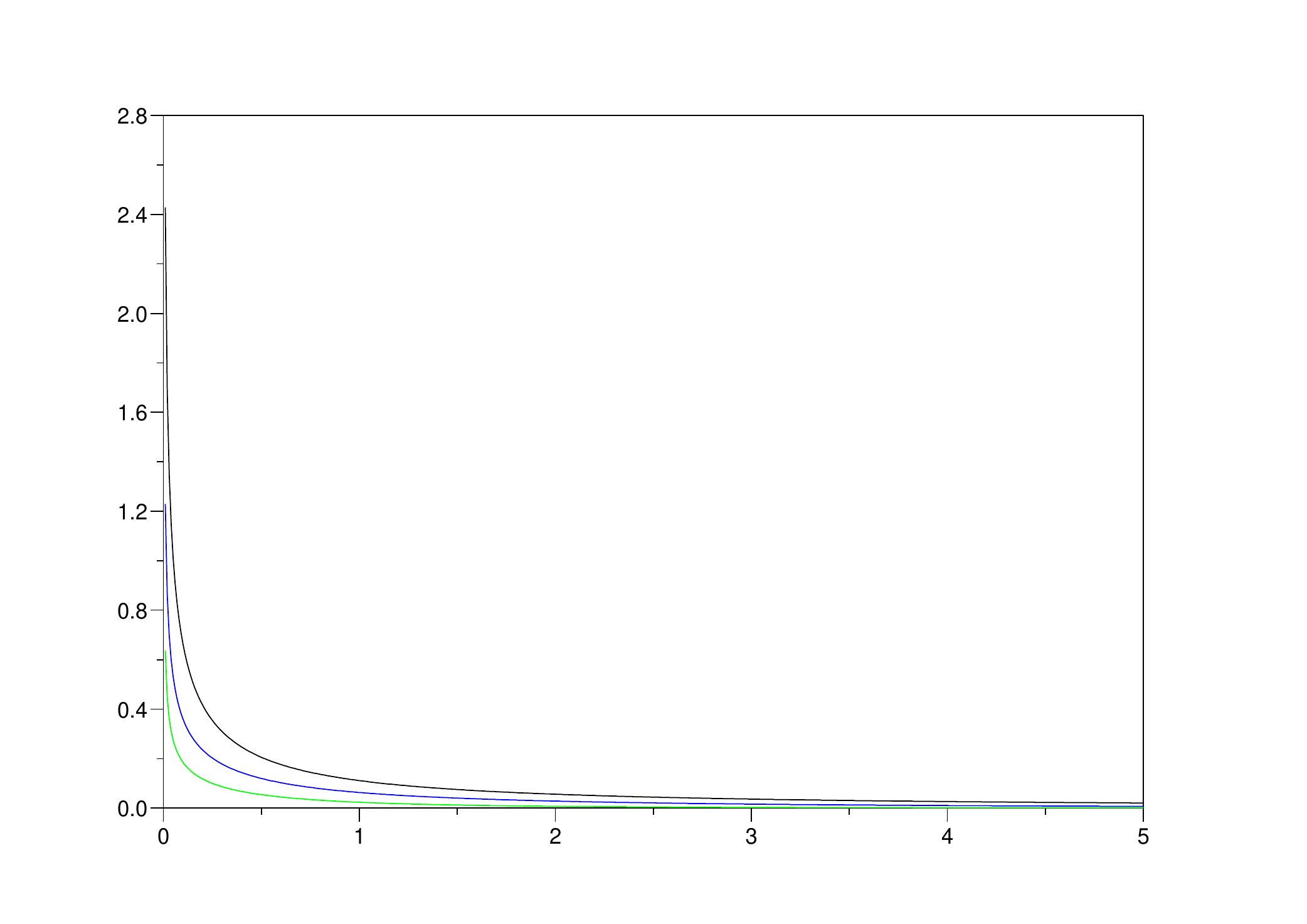}
\caption{Approximate density of $T_2(x)$ for $x=0, \ a =1$ and various values of the parameter $b;$ from top to the bottom: $b=0, \ b=-0.5, \ b=-1.$
}
\end{figure}

\indent As far as the expectation of  $\tau _2 (x)$ is concerned, it is obviously infinite for $b \neq 0 ,$ while $E(\tau _1 (x))$
is finite.
If $b=0, \ E( \tau _1(x))$ and $E( \tau _2(x))$
are both infinite.
As for $\tau _1 (x)$ this is well-known, as for $\tau _2 (x)$ it derives from the fact that (see \cite{abundo:phy16}):
$$ E(T_2 (x)) = \int _ 0 ^{+ \infty } P( T_2(x) >t ) dt = $$
$$ \int _ 0 ^{ + \infty } ds  \frac 2 \pi \ \frac {|a-x| } {\sqrt { 2 \pi } s^{3/2} } \ e^{ - (a -x)^2 /2s } \int _ 0 ^{ + \infty }
\arcsin \sqrt { \frac s {s+t}} dt  = + \infty  .$$
By taking the derivative with respect to $t$ in \eqref{conditionaldistrT2},
one obtains the  density of $T_2(x)$ conditional to $\tau _1(x)=s,$ that is:
\begin{equation} \label{denT2|tau1-BM}
f_{T_2 | \tau _1} (t | s) = - \frac d {dt} \int _0 ^s  \psi _ {s +t } (y) dy = e^ {- b^2 (s+t)/2} \frac {\sqrt s } { \pi (s+t) \sqrt t } \ ,
\end{equation}
and, for $b=0:$
\begin{equation}
f_{T_2 | \tau _1} (t | s) = \frac {\sqrt s } { \pi (s+t) \sqrt t } .
\end{equation}
Since $ \tau _2 (x) = \tau _1 (x) + T_2 (x),$ by the convolution formula, the density of $\tau _2 (x)$ follows:
\begin{equation} \label{dentau2}
f_ {\tau _2} (t) = \int _ 0 ^t f_{T_2 | \tau _1} (t-s | s) f _ {\tau _1 } (s) ds  = \frac {e^ { - b^2 t/2}} { \pi t}
\int _0 ^t \frac {|a -x| } { \sqrt { 2 \pi } s \sqrt { t-s } } e^{ - (a +bs-x)^2 /2s } ds .
\end{equation}
Of course, the distribution of $\tau _2 (x)$  is defective for $b \neq 0,$ namely
$\int _ 0 ^ {+ \infty} f_ {\tau _2} (t) dt = 1 - P( \tau _2 (x)= + \infty) < 1 ,$
since $P( \tau _2 (x)= + \infty) = P( T _2 (x)= + \infty) >0.$ \par\noindent
If $b=0,$ we obtain:
\begin{equation} \label{dentau2b0}
f_ {\tau _2} (t) = \frac 1 { \pi t }
\int _0 ^t \frac {|a -x| } { \sqrt { 2 \pi } s \sqrt { t-s } } e^{ - (a -x)^2 /2s } ds ,
\end{equation}
which is non-defective. \par
In the Figure 3, we report the probability density of $\tau _2 (x),$ obtained from \eqref{dentau2} by calculating numerically the integral, for
$x=0, \ a=1,$ and various values of the parameter $b \le 0.$ Although the shapes appear to be similar to that of the inverse Gaussian density
\eqref{IGdensity}, the density
of $\tau _2 (x)$ is more concentrated around its maximum. In the Figure 4, we report from \cite{abundo:phy16} the comparison between the probability density of $\tau _2 (x)$ and the inverse Gaussian density, for $a=1, b=0$ and $x=0.$

\begin{figure}
\centering
\includegraphics[height=0.33 \textheight]{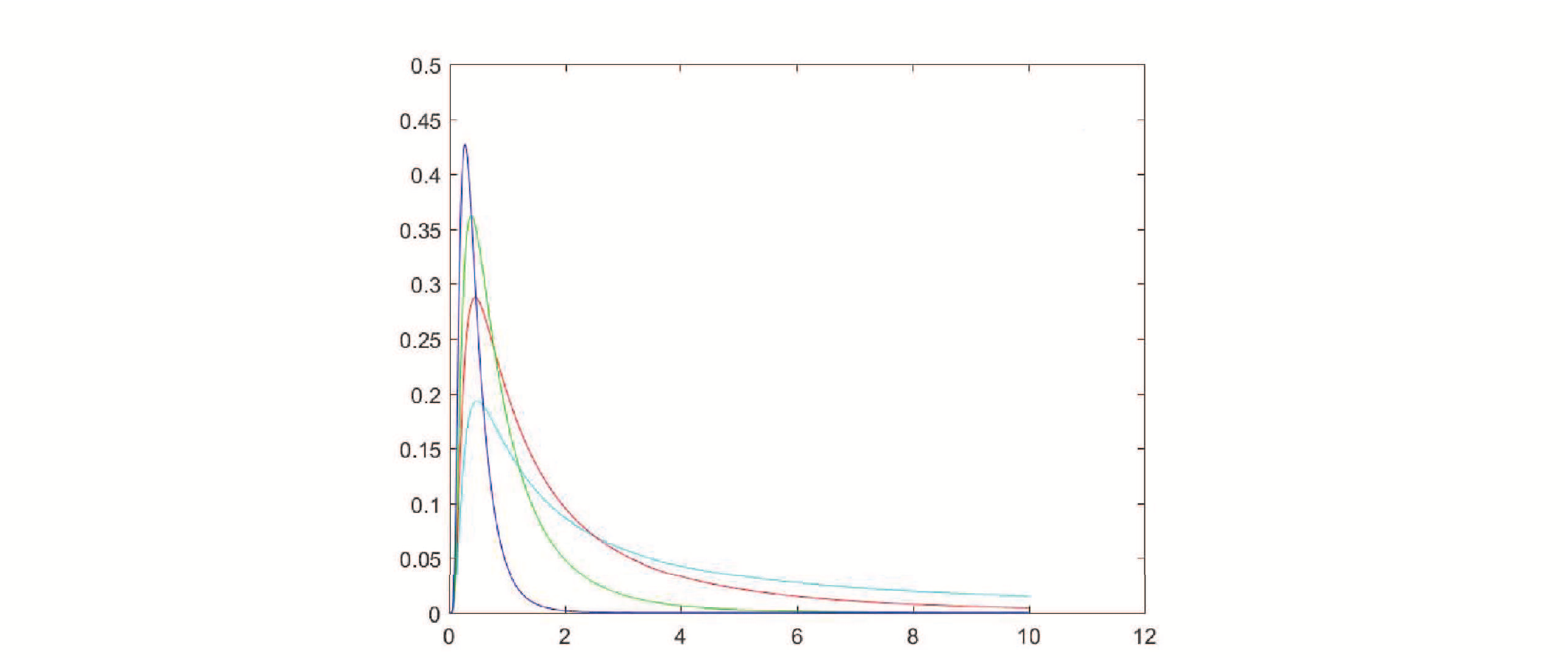}
\caption{Approximate density of $\tau _2(x)$ for $x=0, \ a =1,$ and various values of the parameter $b;$ from top to the bottom,
with respect to the peak of the curve: $b=-2, \ b=-1, \ b=-0.5, \ b=0.$
}
\end{figure}

\begin{figure}
\centering
\includegraphics[height=0.33 \textheight]{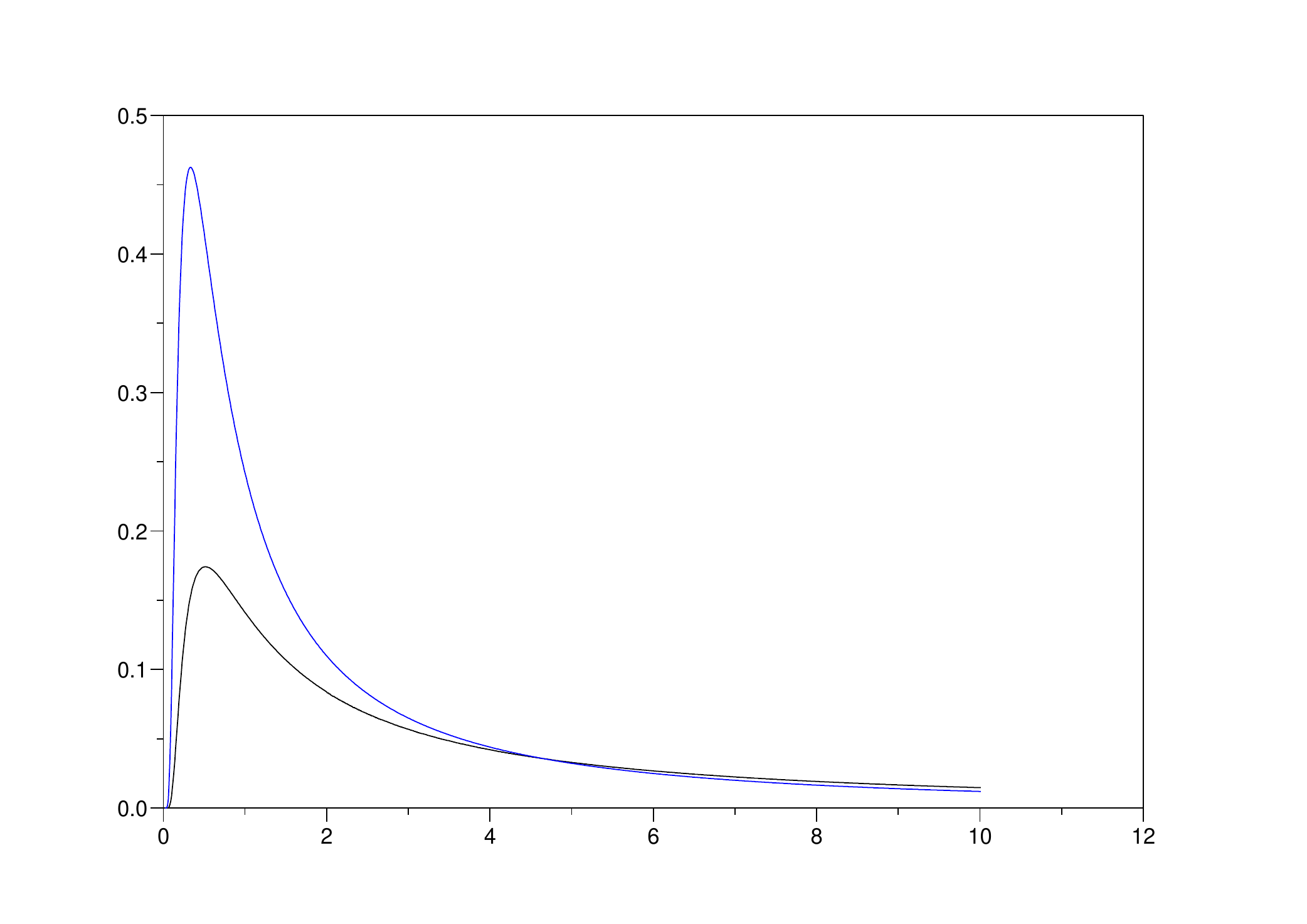}
\caption{Comparison between the probability density of $\tau _2 (x)$ (upper peak) and the inverse Gaussian density (lower peak),
for $a=1, b=0$ and $x=0.$
}
\end{figure}
The following holds (see \cite{abundo:phy16}):
\begin{Proposition} \label{proposizioneuno}
Let be
$T_1(x) = \tau _1 (x), \ T_n (x)= \tau  _n (x) - \tau _ {n-1} (x), \  n =2, \dots $  \par
Then:
\begin{equation}
P(T_1 (x) \le t ) = 2(1 - \Phi(a -x / \sqrt t  )),
\end{equation}
\begin{equation}
P( T_n (x) \le t ) = 1 - \int _ 0 ^{+ \infty } f_{\tau _ {n-1} } (s) ds \int _ 0 ^s \psi _{s+t} (y) dy, \ n= 2, \dots
\end{equation}
Moreover, the density of $\tau _n (x)$ is:
\begin{equation}
f_ {\tau _n} (t) = \int _ 0 ^t f_{T_n | \tau_{n-1}} (t-s | s) f _ {\tau _{n-1} } (s) ds ,
\end{equation}
where $f_{\tau _ {n-1} }$ and $f_ {T_n | \tau _ {n-1} }$ can be calculated inductively, in a similar way, as done for $f_{\tau _ {2} }$
and $f_{T_2 | \tau _1 } .$ \par\noindent
If $b=0, \ T_1(x), \ T_2 (x), ....$ are finite with probability one.
\end{Proposition}

\hfill $\Box$
\bigskip

\noindent We conclude this section, by summarizing the Laplace transforms  of the first and second passage-time of $B_t$ through the boundary $S(t)= a+ bt:$ \par\noindent
$\bullet $ the Laplace transform (LT) of $\tau _1 (x)$ is the LT of the inverse Gaussian density:
\begin{equation}
\psi(\lambda) = E \left [ e^ { - \lambda \tau_1(x) } \right ] =  e^{-(a - x)(\sqrt{b^2 + 2\lambda} - b)}, \ \lambda >0 ;
\end{equation}
$\bullet $ the Laplace transform (LT) of $\tau _2 (x)$ is:
\begin{equation}
\tilde{\psi}(\lambda) = E \left[e^{- \lambda \tau_2(x)}\right] = \int_0^{+\infty}f_{\tau_2}(t)e^{-\lambda t}\;dt =
\int_0^{+\infty}e^{-\lambda t}\left[\frac{e^{-\frac{b^2}{2}t}}{\pi t} \int_0^t \frac{|a - x|e^{-\frac{(a + bs - x)^2}{2s}}}{\sqrt{2\pi}\cdot
s \sqrt{t - s}}\;ds\right]\;dt,
\end{equation}

\begin{figure}
\centering
\includegraphics[height=0.33 \textheight]{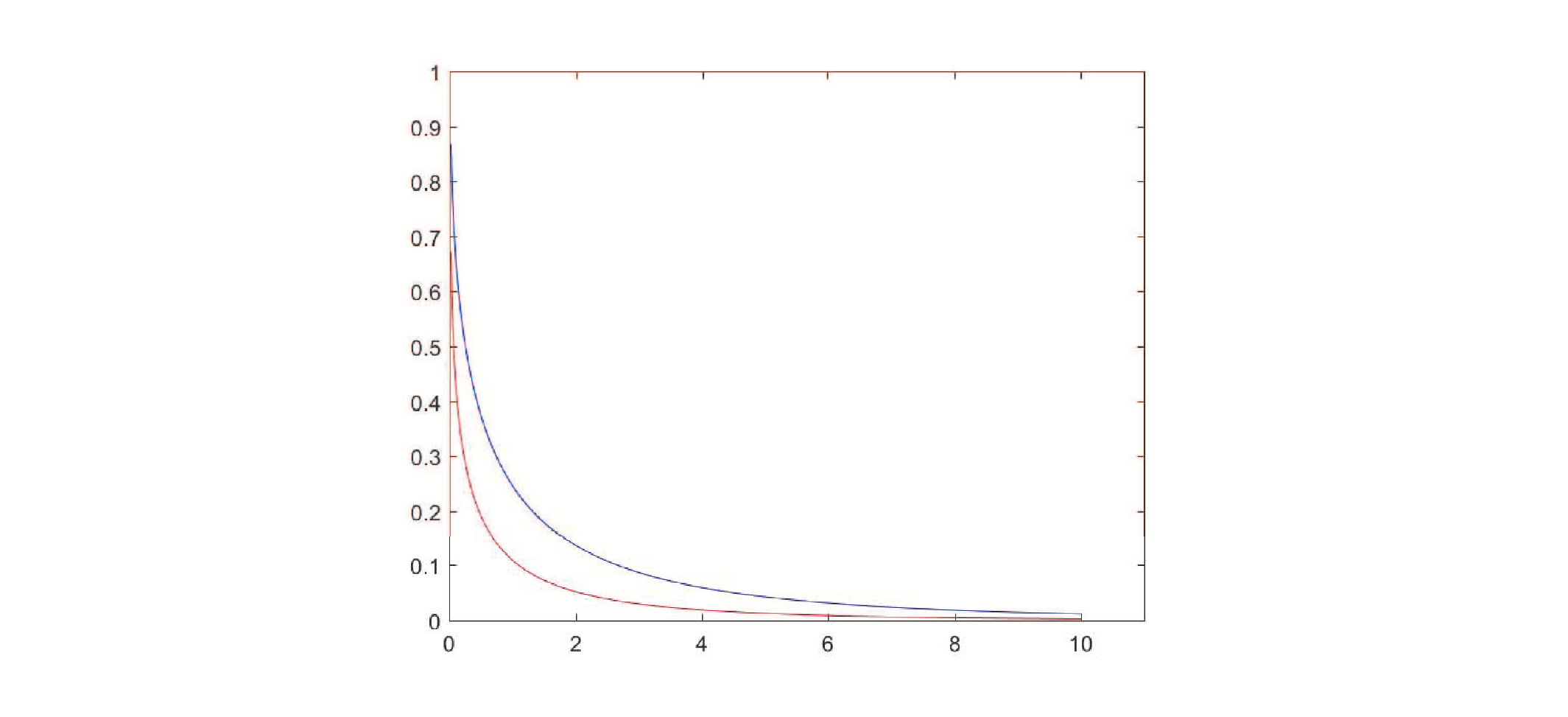}
\caption{Comparison between the  LT of the inverse Gaussian density (upper curve) and the LT  of the second-passage time  $\tau_2(x)$
of BM through the boundary $S(t)=a + bt $ (lower curve),
for $x = 0$ and $a =1, \ b=0$. \label{Fig : 13}
}
\end{figure}

In the Figure \ref{Fig : 13}, we compare the shape of the LT of the inverse Gaussian density with the LT  of the second-passage time  $\tau_2(x)$
of BM through the boundary $S(t)=a + bt ,$ for $a=1, b=0$ and $x=0.$

\section{The nth-passage time of a transformed Brownian motion}
In this section we will extend the previous results  to a one-dimensional, time homogeneous diffusion $X(t),$ obtained from BM by a space
transformation or
a time change. We start with diffusions conjugated to BM.

\subsection{The nth-passage time of a diffusion conjugated to BM}
Let be $X(t)$ a  one-dimensional, time-homogeneous diffusion process, which is driven by the SDE
$$ dX(t) = \mu (X(t)) dt + \sigma (X(t)) d B_t, \ X(0)= x ,$$
where the coefficients $\mu ( \cdot )$ and $ \sigma ( \cdot)$ satisfy the usual conditions for the existence and uniqueness of
the solution (see e.g. \cite{ikeda:81}).\par
We recall that $X(t)$ is said to be {\it conjugated  to BM} (see \cite{abundo:stapro12}), if there exists an increasing function
$v$ with $v(0) = 0$, such that:
$$X(t) = v^{-1}\left(B_t + v(x)\right), \  t \geq 0.$$
Examples of diffusions conjugated to BM are (see also \cite{abundo:stapro12}):
\par\noindent
(i) the Feller process or Cox-Ingersoll-Ross (CIR) model, which is driven by the SDE
$$dX(t) = \frac 14 dt + \sqrt{X(t) \lor 0}\;dB_t, \qquad X(0) = x,$$ and is
conjugated to BM via the function $v(x) = 2\sqrt{x},$ i.e.
$X(t) = \frac 14\left(B_t + 2\sqrt{x}\right)^2;$ \par\noindent
(ii) the Wright and Fisher-like process, which is driven by the SDE
$$dX(t) = \left(\frac 14 - \frac 12 X(t)\right)\;dt + \sqrt{X(t)\left(1 - X(t)\right) \lor 0} \; dB_t, \qquad X(0) = x \in [0, 1],$$
and is conjugated to BM via the function $v(x) = 2 \; \arcsin\sqrt{x}$, i.e.
$$X(t) = \sin^2\left(B_t/2 + \arcsin\sqrt{x}\;\right).$$ \par
Suppose now that $X(t)$ is conjugated to BM via the function $v;$
let $\tau _1 (a|x) = \inf \{ t>0: X(t) =a \}$ be
the FPT of $X(t)$ through
the level $ a >x ;$
as easily seen, one has
$ \tau _1 (a|x) =  \tau _1^{B}(v(a) |v(x)),$
where $ \tau _1^{B}(\alpha |y),$ denotes the FPT of BM, starting from $y < \alpha ,$ through the level $\alpha$.
Then, from \eqref{IGdensity}, one gets the density of $\tau _1 (a|x):$
\begin{equation}
f_{\tau_1}(t) =  \frac{v(a) - v(x)}{t^{3/2}} \phi\left(\frac{v(x) - v(a)}{\sqrt{t}}\right),
\end{equation}
implying that
$ E \left[ \tau _1 (a|x) \right] = +\infty.$
For $x < a,$ one has:
$$ \tau _2 (a|x)  = \inf\{t > \tau _1 (a|x): X(t) = a \} =
\inf\{t > \tau _1 (a|x): B_t = v(a) - v(x) \} =  \tau _2^B (v(a)|v(x)) ,$$
and
$$\tau _n (a|x) = \inf\{t > \tau_{n - 1}(a|x) : B_t = v(a) - v(x) \} = \tau _n^B (v(a)|v(x)) ,$$
where the superscript $B$ refers to BM, namely $\tau_n^B(a|x), \ n=2, \dots $ is the $n$th-passage time of BM through the level $a,$
when starting from $x <a .$
By  using the same notation of the previous section, setting $T_1(x) = \tau_1 (a |x) = \tau _1 ^B (v(a)|v(x)), \  T_2(x) = \tau _2 (a|x) - \tau _1 (a|x) =
\tau _2 ^B (v(a)|v(x)) - \tau _1 ^B (v(a)|v(x)), $
and
using \eqref{denT2} with $b=0,$ we get
\begin{equation}
P\left(T_2(x) \leq t\right) = \frac
2\pi\int_0^{+\infty}f_{\tau_1}(s) \arccos\sqrt{\frac s{s + t}}\;ds.
\end{equation}
Thus, the density of $T_2$ is
$$
f_{T_2}(t)= \int_0^{+\infty}
f_{\tau_1}(s) \ \frac{\sqrt{s}}{\pi\sqrt{t} (s + t)}\;ds  = \int_0^{+\infty}\frac{v(a) - v(x)}{s^{3/2}} \
\frac{e^{-\frac{(v(x) -
v(a))^2}{2s}}}{\sqrt{2\pi}} \ \frac{\sqrt{s}}{\pi\sqrt{t}(s + t)}\;ds $$
\begin{equation} \label{fT2}
= \frac{v(a) - v(x)}{\pi\sqrt{2\pi}}\int_0^{+\infty} \frac{1}{(s +
t)\sqrt{t}} \ \frac 1s \  e^{-\frac{(v(x) - v(a))^2}{2s}}\;ds.
\end{equation}
The expectation of $\tau_2 (a|x)$ is infinite,
because $E \left [T_2(x) \right] = +\infty .$
From \eqref{denT2|tau1-BM} with $b=0,$ we get
\begin{equation}
f_{T_2 \arrowvert \tau_1}(t\arrowvert s) = -\frac{d}{dt}\int_0^s\psi_{s + t}(y)\;dy = \frac{\sqrt{s}}{\pi(s + t)\sqrt{t}}.
\end{equation}
Moreover, from \eqref{dentau2b0} $(b=0):$
\begin{equation}
f_{\tau_2}(t) = \int_0^t\frac{v(a) - v(x)}{\pi ts\sqrt{t -
s}\sqrt{2\pi}}\cdot e^{-\frac{(v(x) - v(a))^2}{2s}}\;ds.
\end{equation}
Setting $T_n (x)= \tau_n(a|x) - \tau_{n-1} (a|x),$ and
using Proposition \ref{proposizioneuno}, we get
\begin{equation}
P \left(T_1(x) \leq t\right) = 2\left(1 - \Phi\left(\frac{v(a) - v(x)}{\sqrt{t}}\right)\right),
\end{equation}
\begin{equation}
P \left(T_n(x) \leq t \right) = 1 -
\int_0^{+\infty}f_{\tau_{n - 1}}(s)\;ds\int_0^s\psi_{s + t}(y)\;dy, \quad n = 2, \dots.
\end{equation}
The density of $\tau_n (a|x)$ is:
\begin{equation}
f_{\tau_n}(t) = \int_0^t f_{T_n | \tau_{n - 1}}(t - s|  s)\cdot f_{\tau_{n - 1}}(s)\;ds,
\end{equation}
where $f_{\tau_{n - 1}}$ and $f_{T_n |  \tau _{n - 1}}$ can be calculated inductively, in a similar way, as
done for $f_{\tau _2}$ and $f_{T_2 | \tau_1}$. In conclusion, $T_1(x), T_2(x), \dots$ are finite with probability one.

\subsection{The nth-passage time of time-changed Brownian motion}
The previous arguments can be applied  also to time-changed BM, namely
\begin{equation}\label{time-changed}
X(t) = x + B(\rho(t)),
\end{equation}
where $\rho(t) \geq 0$ is an increasing, differentiable function of $t > 0$, with $\rho(0) = 0$.
Such kind of diffusion process $X(t)$ is a special case of Gauss-Markov process; in particular the form (\ref{time-changed})
is taken by certain integrated
Gauss-Markov processes (see \cite{abundo:smj13}, \cite{abundo:smj15}).
Let us consider the constant barrier $S = a,$ and $x <a;$ then, the FPT of $X(t)$ through $a$ is
$ \tau _ 1 (a|x) = \rho^{-1}( \tau_1 ^B (a|x)),$
where $\tau _1 ^B (a|x)$ denotes the first-passage time of BM, starting from $x$, through the barrier $a$.
Thus, the density of $\tau _ 1 (a|x)$ is :
\begin{equation} \label{FPT-TC}
f_{\tau_1}(t) = f_{\tau _1^B }(\rho(t))  \rho'(t) = \frac{|a - x|}{\rho(t)^{3/2}} \phi\left(\frac{a - x}{\sqrt{\rho(t)}}\right) \rho'(t).
\end{equation}
The expectation is:
\begin{equation} \label{mediaTau1_TC}
E \left[\tau _ 1 (a|x)\right] = \int_0^{+\infty}t \  f_{\tau_1}(t)\;dt = \int_0^{+\infty} t \frac{|a - x|\cdot
\rho'(t)}{\rho(t)^{3/2}} \frac{e^{-\frac{(x - a)^2}{2\rho(t)}}}{\sqrt{2\pi}}\;dt
\end{equation}
Suppose that
$ \rho(t) \sim c\cdot t^{\alpha} ,$ as $ t \rightarrow + \infty ;$ then, as easily seen, for $\alpha >2$ the integral converges,
namely $E \left[\tau _ 1 (a|x)\right] < +\infty,$
unlike the case of BM, for which the expectation of FPT is infinite.
One has:
$$
\tau_2(a|x) = \rho^{-1}(\tau _2 ^B (a|x)) \ {\rm with} \; \tau _2 ^B (a|x) = \inf\{s > \tau _1 ^B (a|x) : x+B_s = a \},
$$
$ \ \ \ \ \ \ \ \ \ \ \ \ \ \dots \dots \dots \dots  ,$
$$ \tau_n(a|x) =
\rho^{-1}(\tau _n ^B (a|x)) \ {\rm with} \; \tau _n ^B (a|x) = \inf\{s > \tau _{n-1} ^B (a|x) : X+ B_s = a \};
$$
$$
T_1(x) = \tau_1(a|x) = \rho^{-1}(\tau _1 ^B (a|x)),  \ T_2(x) = \tau_2(a|x) - \tau_1(a|x) =
\rho^{-1}(\tau _2 ^B (a|x)) - \rho^{-1}(\tau _1 ^B (a|x)).
$$
It results
\begin{equation}
P\left(T_2(x) \leq t\right) = \int_0^{+\infty} P\left(\rho^{-1}\left(\tau _2 ^B (a|x)\right) \leq t + s\right)
f_{\rho^{-1}(\tau _ 1^B (a|x))}(s)\;ds ,
\end{equation}
and so the density of $f_{T_2}$ is
\begin{equation} \label{fT2_TC}
 f_{T_2}(t) = \int_0^{+\infty}f_{\tau _2 ^B }\left(\rho(t + s)\right) \rho'(t + s) f_{\tau_1}(s)\;ds
= \int_0^{+\infty}f_{\tau _2 ^B}\left(\rho(t + s)\right) \rho'(t + s) f_{\tau _1 ^B}(\rho(s)) \rho'(s)\;ds.
\end{equation}

Moreover,
$$ P\left(T_2(x) = +\infty\right) = \mathop{\lim}_{t \to +\infty} P \left(T_2(x) > t\right) $$
$$= 1 - \int_0^{+\infty}f_{\tau_1}(s) \left[\mathop{\lim}_{t \to +\infty}P \left(\tau _2 ^B (a|x) \leq \rho(t + s)\right)\right] \ ds .
$$
Therefore, by using \eqref{dentau2b0}:
$$ P\left(T_2(x) = +\infty\right) = 1 - \int_0^{+\infty}f_{\tau_1}(s)
\left[\mathop{\lim}_{t \to +\infty}\int_0^{ \rho (t+s)}f_{\tau_2^B}(u)\;du \right]\;ds
$$
$$ =  1 - \int_0^{+\infty}f_{\tau_1}(s)\bigg[\mathop{\lim}_{t \to +\infty}\int_0^{+\infty}\bigg[\frac 1{\pi u}
 \int_0^{\rho(t + s)}\frac{|a - x| e^{-\frac{(a - x)^2}{2y}}}{\sqrt{2\pi} y \sqrt{u - y}}\;dy
\bigg]\;du\bigg]\;ds . $$
Moreover
$$
E\left[T_2(x)\right] = \int_0^{+\infty} P \left(T_2(x) > t\right)\;dt $$
\begin{equation} \label{MediaT2}
= \int_0^{+\infty}\left[1 -
\int_0^{+\infty} P (\tau_2 ^B(a|x) \leq \rho(t + s)) f_{\tau_1}(s)\;ds\right],
\end{equation}
which, unfortunately, it is impossible to calculate, due to the complexity of the integral.\par
The density of $ \tau _2 (a|x)$ is
$$
f_{\tau_2}(t)  =
\frac d{dt}P\left(\tau_2 ^B (a|x) \leq \rho(t)\right) = f_{\tau_2^B}(\rho(t))\cdot \rho'(t) $$
\begin{equation} \label{fTau2_TC}
= \frac{\rho'(t)}{\pi \rho(t)} \int_0^{\rho(t)}\frac{|a - x|}{\sqrt{2\pi} s \sqrt{\rho(t) - s}} e^{-\frac{(a -
x)^2}{2s}}\;ds.
\end{equation}
By using that
$ T_n(x) = \tau_ n(a|x) - \tau _{n - 1}(a|x) =\rho^{-1}(\tau _n^B (a|x)) - \rho^{-1}(\tau_{n - 1}^B (a|x)),$
and conditioning with respect to $\tau_ {n - 1}(a|x) = s$, one gets:
\begin{equation}
P\left(T_n(x) \leq t\right) = \int_0^{+\infty} P\left(\rho^{-1}(\tau _n ^B(a|x)) \leq t + s\right) f_{\tau_{n - 1}}(s)\;ds .
\end{equation}
Thus, the density of $T_n (x)$ is:
$$
f_{T_n}(t) = \int_0^{+\infty}\frac d{dt} P(\tau _n ^B(a|x) \leq \rho(t + s)) f_{\tau_{n - 1}}(s)\;ds =
\int_0^{+\infty}f_{\tau_n ^B}(\rho(t + s))\rho'(t + s) f_{\tau^{n - 1}}(s)\;ds ,
$$
while the density of $\tau_ n(a|x)$ is :
$$
f_{\tau_ n}(t) = \frac d{dt} P \left(\tau_n ^B (a|x) \leq \rho(t)\right) $$
$$=
f_{\tau_n^B}(\rho(t))\rho'(t) = \int_0^{\rho(t)}f_{T_n^B | \tau_{n - 1,^B}}(\rho(t) - s | s) f_{\tau_ {n - 1}^B}(s)\;ds \cdot
\rho'(t),
$$
where $f_{T_n^B | \tau_{n - 1} ^B}$ and $f_{\tau_{n - 1} ^B}$ can be calculated inductively.

\section{Some examples}
Now, we will study  two interesting cases of time-changed BM, with regard to successive-passage times: $B(\frac{t^3}3)$ and the Ornstein-Uhlenbeck (OU) process.
\bigskip

\noindent {\bf 1) \ $B( t^3 /3)$ and integrated BM}. \par\noindent
Let be
\begin{equation}
X(t) = \int_0^t B_s\;ds
\end{equation}
the integrated BM, which has many important application (see e.g. \cite{abundo:smj13}, \cite{abundo:smj15}). It is a Gaussian process with mean zero and covariance function \par\noindent
$Cov(X(t), X(s)) = s^2 (t/2 - s/6), \ s \le t $ (see e.g. \cite{ross:pmodel10}).
Thus, its variance is $Var (X(t)) = t^3/3,$
and $X(t)$ has the same distribution as $B(\frac{t^3}3)$.
Motivated from this, we study the successive-passage time of $B(\frac{t^3}3)$ through a constant barrier $a$.
Since the starting point is zero, we consider
$$\tau _1 (a|0)= \inf\{t > 0 : X(t) = a \} = \inf \{t > 0 : B (t^3 / 3  ) = a\},$$
$$\tau _2(a|0) = \inf\{t > \tau_ 1(a|0) : B  (t^3/ 3  ) = a \},$$
and, as before,
$$T_1(x) = \tau_1(a|0), \ T_2(x) = \tau_ 2(a|0) - \tau_ 1(a|0).$$
Then, by using \eqref{FPT-TC}, we get the density of $\tau_1 (a|0):$
\begin{equation} \label{fTau1_BI}
f_{\tau_1}(t) = \frac{|a - x|}{\left(\frac{t^3}3\right)^{3/2}}\cdot\phi\left(\frac{a - x}{\sqrt{\frac{t^3}{3}}}\right) \ t^2 =
\frac{3\sqrt{3} \ |a - x|}{t^{5/2}} \ \frac{e^{-\frac{3(a - x)^2}{2t^3}}}{\sqrt{2\pi}} .
\end{equation}

\begin{figure}
\centering
\includegraphics[height=0.33 \textheight]{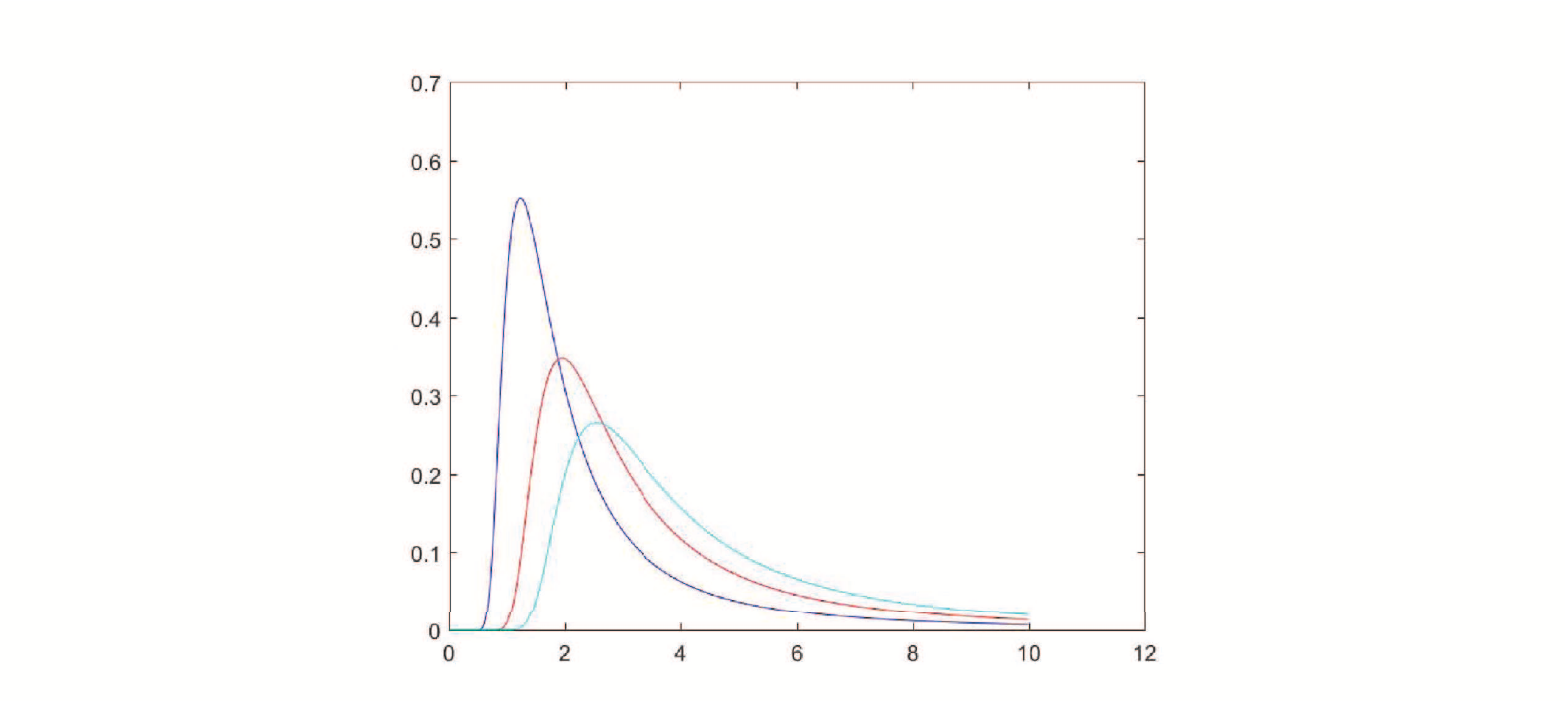}
\caption{Comparison between the shapes of the FPT density of $B(t^3/3)$ \eqref{fTau1_BI}, for different values of the barrier $a$;
from top to the bottom, with respect to the peak of the curve: $a = 1, a = 2, a = 3.$ \label{Fig : 5}
}
\end{figure}

In the Figure \ref{Fig : 5} we report the FPT density of $B(\frac{t^3}3)$ obtained from \eqref{fTau1_BI},
for various values of the barrier: $a = 1, \ a = 2,  \ a = 3$.
The expectation of $\tau_ 1(a|0)$ is obtained by \eqref{mediaTau1_TC}:
\begin{equation}
E\left[\tau_1(a|0)\right] = \frac{3\sqrt{3} \ |a - x|}{\sqrt{2\pi}} \int_0^{+\infty}\frac{e^{-\frac{3(a - x)^2}{2s}}}{t^{3/2}}\;dt .
\end{equation}
Notice that, unlike the case of BM, it is finite (in fact $\rho(t) = t^3/3
\sim c\cdot t^{\alpha} \mbox{ with } \alpha = 3 > 2$ (see the observation after equation \eqref{mediaTau1_TC}).
Calculating numerically the integral e.g.
for $a = 1$, we have obtained  $E\left[\tau_1(a|0)\right] \approx  3.594 \ ,$ which is
not far from the estimate $\hat{E}\left[\tau_1(a|0)\right] = 3.704 \ ,$ obtained by Monte Carlo simulation.
By using \eqref{fT2_TC}, one has:
\begin{align}
f_{T_2}(t) &= \int_0^{+\infty}\frac{3|a - x|}{\pi (t + s)^{3}} \left[\int_0^{\frac{(t + s)^3}3}\frac{e^{-\frac{(a - x)^2}{2u}}}{u\sqrt{\frac{(t + s)^3}3
- u}}du\right]{(t + s)^2}s^2\frac{3\sqrt{3}|a - x|}{s^{5/2}}\frac{e^{-\frac{3(a - x)^2}{2s^3}}}{\sqrt{2\pi}}ds\notag \\ &= \frac{9\sqrt{3}(a -
x)^2}{2\pi^2}\int_0^{+\infty}\frac{e^{-\frac{3(a - x)^2}{2s^3}}}{(t + s)\sqrt{s}}\left[\int_0^{\frac{(t + s)^3}3}\frac{e^{-\frac{(a - x)^2}{2u}}}{u\sqrt{\frac{(s
+ t)^3}3 - u}}du\right]ds.\label{fT2_BI}
\end{align}

\begin{figure}
\centering
\includegraphics[height=0.33 \textheight]{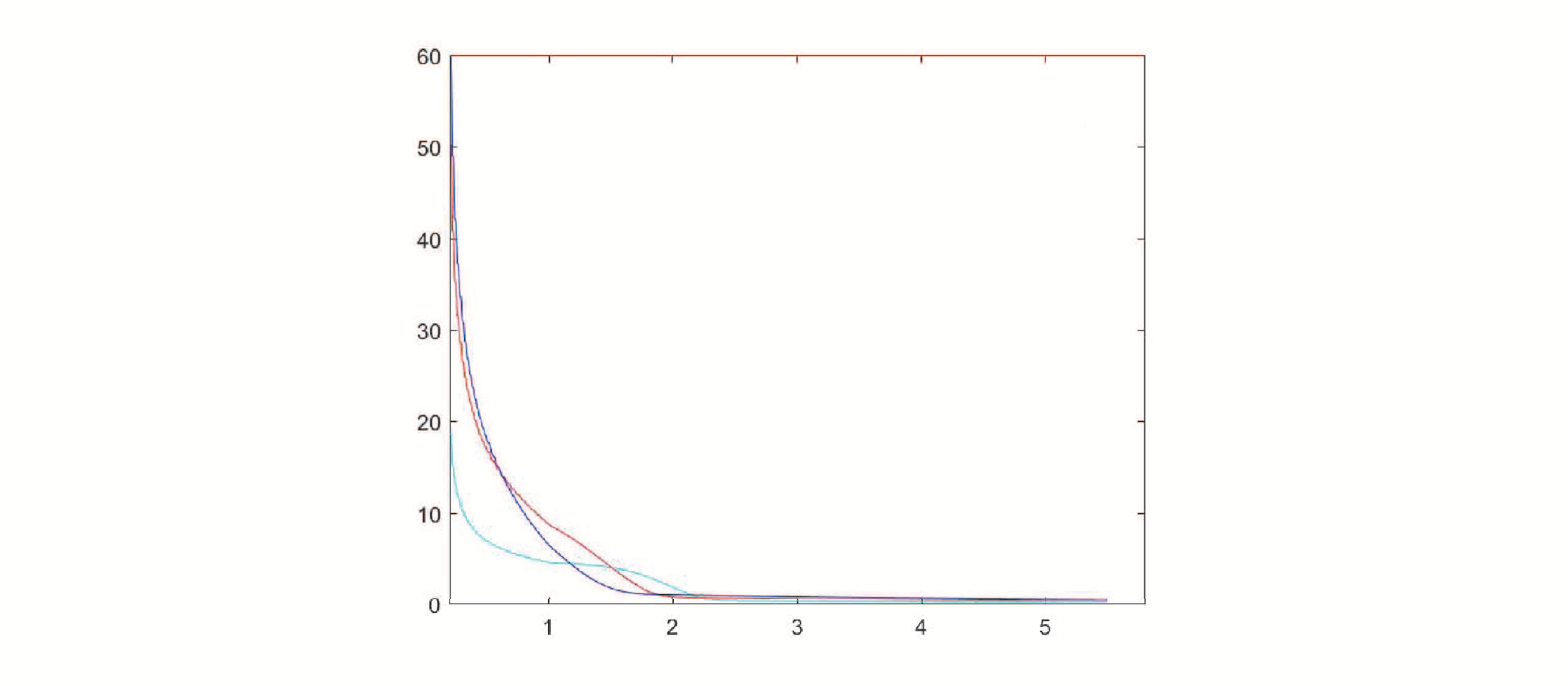}
\caption{Approximate density of $T_2(x)$ \eqref{fT2_BI}, in the case of $B(t^3/3),$ for various values of the parameter $a;$ from top to
bottom, with respect to the value at $t=0:$ $a = 3, \ a = 2, \ a = 1.$ \label{Fig : 6}
}
\end{figure}

In the Figure \ref{Fig : 6}, we report the probability density of $T_2(x)$, obtained from \eqref{fT2_BI} by numerical computation, for
various values of the parameter $a.$

\begin{figure}
\centering
\includegraphics[height=0.33 \textheight]{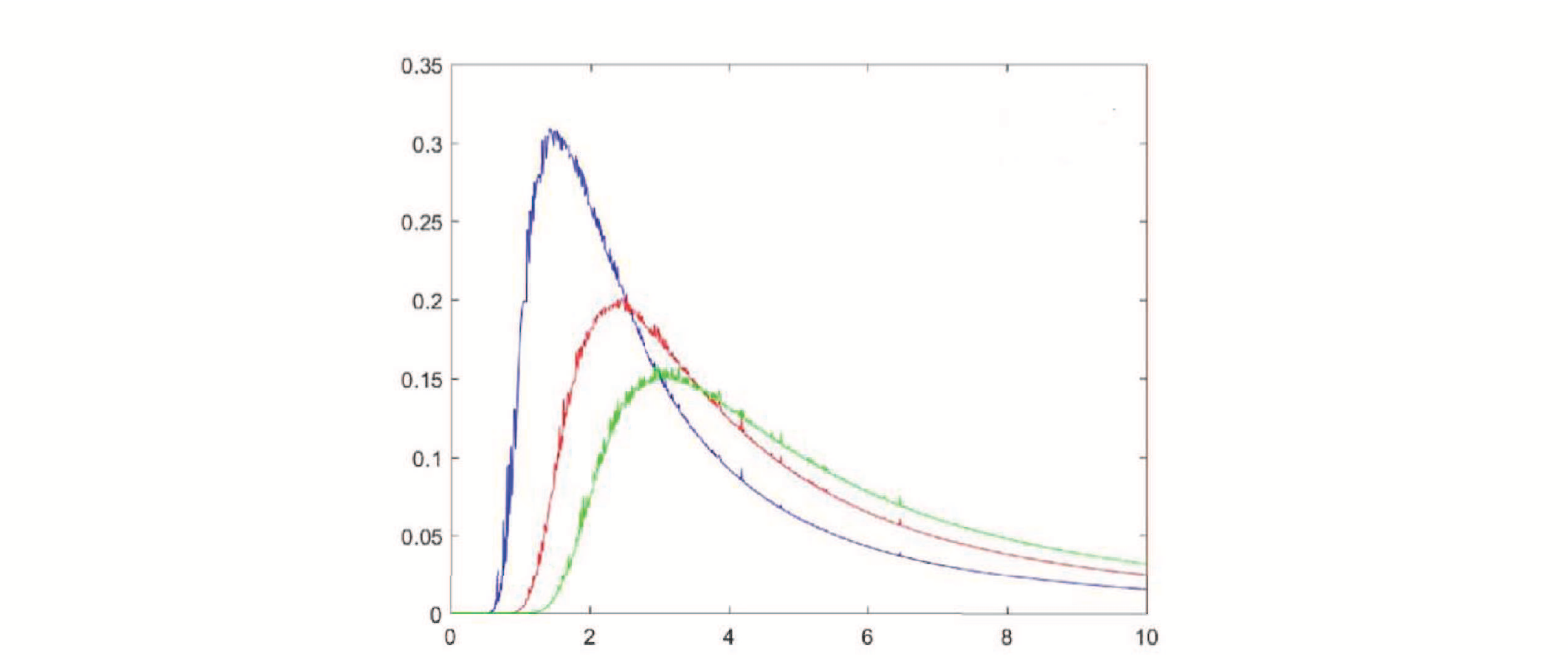}
\caption{Approximate density of $\tau_2(a|x)$ in the case of $B(t^3/3),$ for various values of the parameter $a$; from top to the bottom: $a =
1, a = 2, a = 3:$ \label{Fig : 7}
}
\end{figure}

By using \eqref{fTau2_TC}, one gets the density of $\tau_2 (a|x):$
\begin{equation}\label{fTau2_BI}
f_{\tau_ 2}(t) = \frac 3{\pi t} \int_0^{\frac{t^3}3}\frac{|a - x|}{\sqrt{2\pi}  s \sqrt{\frac{t^3}3 - s}} e^{-\frac{(a - x)^2}{2s}}\;ds .
\end{equation}
Calculating numerically the integral, we obtained the shape  of the density of $\tau_ 2(a|x) ,$ for different values of the parameter $a;$ the results are shown in the Figure \ref{Fig : 7}.

As for the expectation of $\tau_2(a|0),$ one has
$E\left[\tau_ 2(a|0)\right] = - \frac {d } {d \lambda } \tilde{\psi ( \lambda )} \big{\arrowvert}_{\lambda = 0},$
where  $\tilde{\psi}(\lambda)$ is  the Laplace transform of $\tau_2(a|x)$, i.e.
\begin{equation}
\tilde{\psi}(\lambda) =  \int _0 ^ { + \infty } e^ { - \lambda t } f_ { \tau _2} (t) dt = \frac{3|a - x|}{\pi\sqrt{2\pi}} \int_0^{+\infty}\frac{e^{-\lambda t}}{t} \left[\int_0^{\frac{t^3}3}\frac{e^{-\frac{(a -
x)^2}{2s}}}{s\sqrt{\frac{t^3}3 - s}}ds\right]dt.
\end{equation}
Thus:
\begin{equation} \label{MediaTau2_BI}
E\left[\tau_ 2(a|x)\right] = \frac{3|a - x|}{\pi\sqrt{2\pi}} \int_0^{+\infty}\left[\int_0^{\frac{t^3}3} \frac{e^{-\frac{(a -
x)^2}{2s}}}{s\sqrt{\frac{t^3}3 - s}}\;ds\right]\;dt .
\end{equation}
For not too large values of $a$, Monte Carlo  simulation provides the estimate $\widehat{E}\left[\tau_2(a|0)\right]$, not  far from the
exact value, obtained by using \eqref{MediaTau2_BI}.

\begin{figure}
\centering
\includegraphics[height=0.33 \textheight]{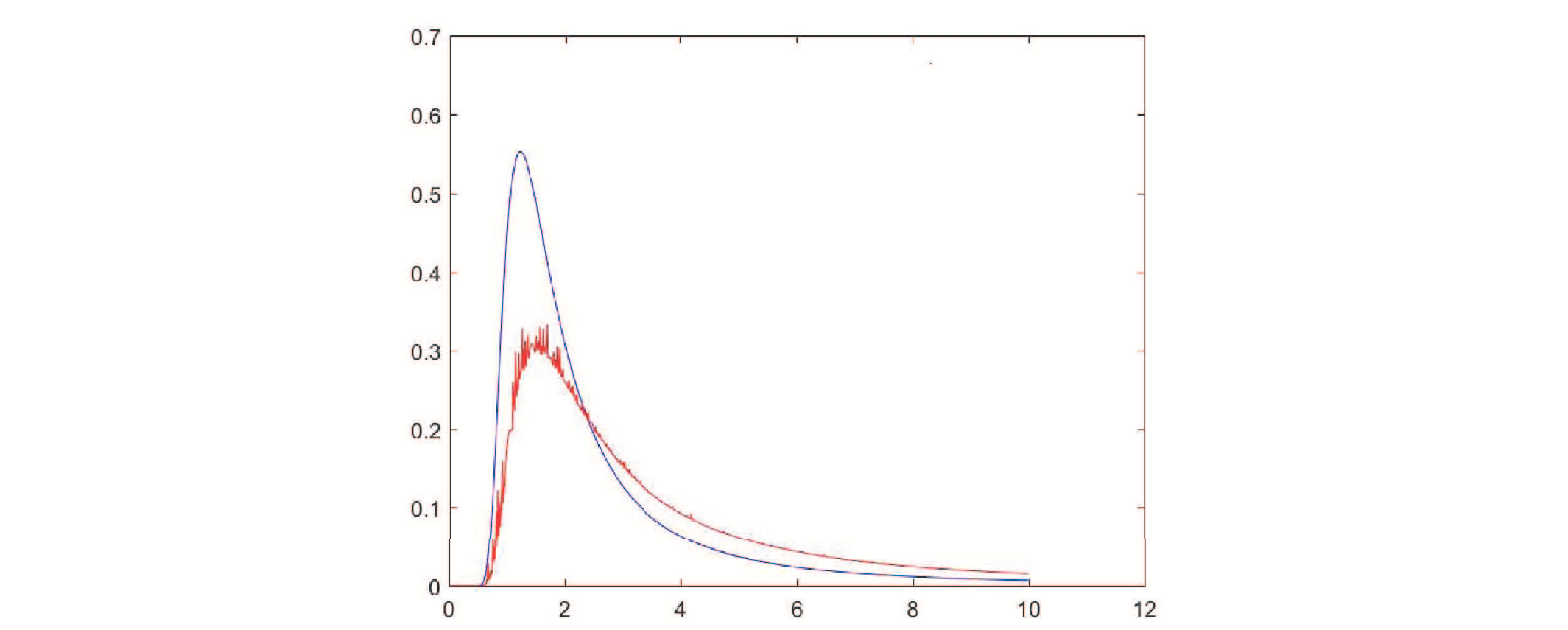}
\caption{Comparison between the density of $\tau_2 (a|0)$ (lower peak) and the
density of $\tau_1(a|0)$ (upper peak), in the case of $B(t^3/3),$ for $a = 1$. \label{Fig : 8}
}
\end{figure}

In the Figure \ref{Fig : 8}, we report the comparison between the density, \eqref{fTau2_BI}, of $\tau_2(a|0)$ and the density, \eqref{fTau1_BI},
of $\tau_ 1(a|0)$, for  $a = 1 .$
Notice that in both the Figures \ref{Fig : 7} and \ref{Fig : 8} the curves appear to be somewhat jagged, due to approximation in
numerical computation. This also happens in the next Figure 13 and Figure 15.
\newpage

\noindent We conclude the example of $B(t^3/3),$  reporting the LTs of $\tau_1(a|0)$ and $\tau _2 (a|0).$ \bigskip

\noindent $\bullet $ LT of $\tau_1(a|x):$

\begin{equation}
\psi(\lambda) = E\left[e^{- \lambda \tau_1(a|x) }\right] = \int_0^{+\infty}f_{\tau_1}(t)e^{-\lambda t}\;dt =
\int_0^{+\infty}\frac{3\sqrt{3}|a - x|}{t^{5/2}}  e^{-\lambda t}\cdot\frac{e^{-\frac{3(a - x)^2}{2t^3}}}{\sqrt{2\pi}}\;dt;
\end{equation}
$\bullet $ LT of $\tau_2(a|x):$
\begin{equation}
\tilde{\psi}(\lambda) = E\left[e^{- \lambda \tau_2(a|x) }\right] = \int_0^{+\infty}f_{\tau_2}(t)e^{-\lambda t}\;dt =\frac{3|a -
x|}{\pi\sqrt{2\pi}} \int_0^{+\infty}\frac{e^{-\lambda t}}{t} \left[\int_0^{\frac{t^3}{3}}\frac{e^{-\frac{(a - x)^2}{2s}}}{s\sqrt{\frac{t^3}3 -
s}}\;ds\right]\;dt .
\end{equation}

\begin{figure}
\centering
\includegraphics[height=0.33 \textheight]{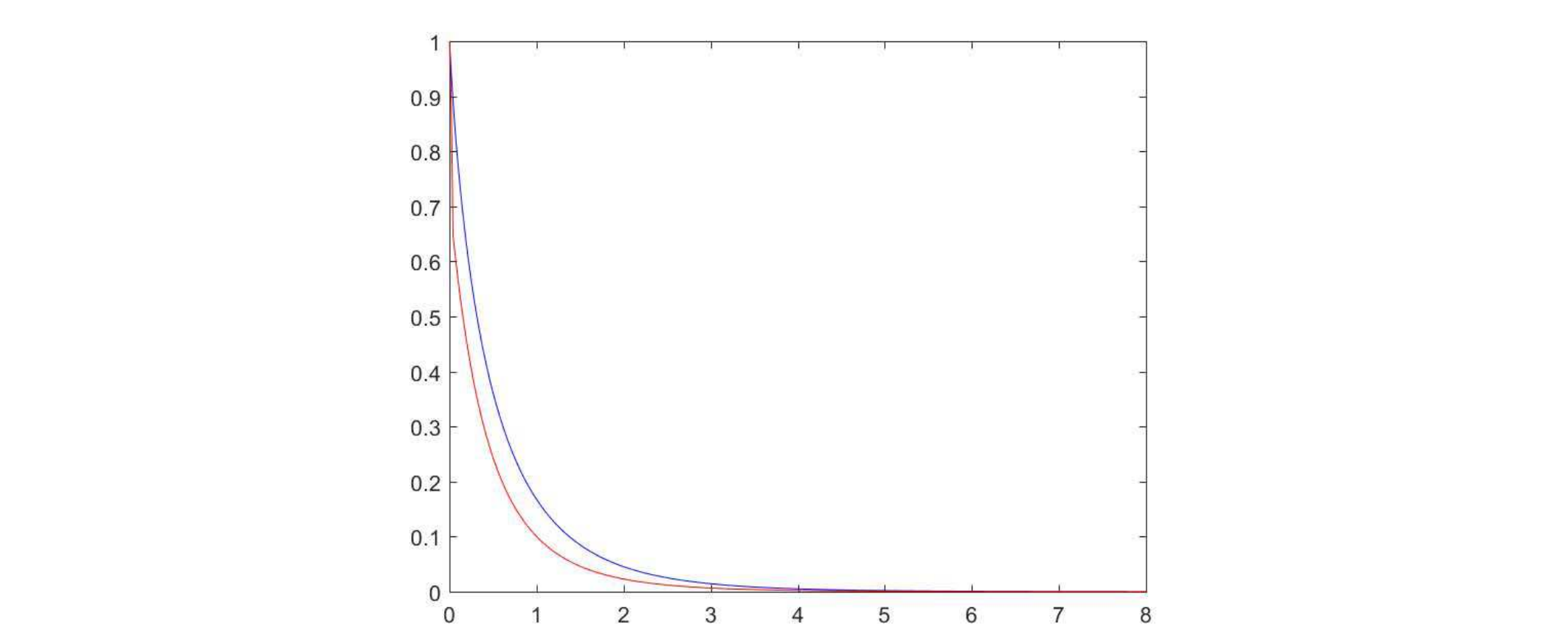}
\caption{Comparison between LT of $\tau_1(a|x)$  (upper curve) and LT of $\tau _2(a|x)$ (lower curve), in the case of
$B(\frac{t^3}3)$, for $a=1$. \label{Fig : 14}
}
\end{figure}

\noindent {\bf 2) The Ornstein-Uhlenbeck (OU) process} \par\noindent
It is the solution of:
\begin{equation}\label{OU}
\begin{cases}
dX(t) = -\mu X(t)dt + \sigma dB(t) \\ X(0) = x
\end{cases}
\end{equation}
with $\mu, \ \sigma $ positive constants. Explicitly, one has
\begin{equation}\label{OUsolution}
X(t) = e^{-\mu t}\left(x + \int_0^t \sigma e^{\mu s}\;dB(s)\right).
\end{equation}
Moreover, (see e.g. \cite{abundo:stapro12}) $X(t)$ can be written as $e^ { - \mu t } \left ( x + B( \rho (t)) \right ),$
where $ \rho (t)= \frac {\sigma ^2 } { 2 \mu}  (1 - e^ {-2 \mu t } ),$ namely in terms of time-changed BM.
Therefore, $X(t)$ has normal distribution with mean $xe^{-\mu t}$ and variance $\frac{\sigma^2}{2\mu}(1 - e^{-2\mu t}) .$
\par\noindent
We choose the time-varying barrier $S(t) = a  e^{-\mu t}$, with $x < a $, and  we  reduce to the passage-times of BM through the constant
barrier $a$. In fact:
\begin{equation}
\tau_1(S|x) = \inf\{t > 0 : X(t) = S(t)\} =\inf\{t > 0 :
x+ B(\rho(t)) = a \}
\end{equation}
and
\begin{equation}
\tau_1(S|x) = \rho^{-1}\left( \tau_1 ^B (a|x)\right),
\end{equation}

where $\tau_1 ^B (a|x)$ is the first passage time of BM, through the barrier $a.$
Moreover,
$$
\tau_2(S|x) = \rho^{-1}\left(\tau_2 ^B (a|x)\right), \dots , \tau_n(S|x) = \rho^{-1}\left(\tau_n ^B (a|x) \right),
$$
which are the successive-passage times of BM through $a.$
Thus, by using \eqref{FPT-TC}, we get the density of $\tau _1 (S|x):$
\begin{equation}\label{fTau1_OU}
f_{\tau_ 1}(t) = \frac{2\mu\sqrt{2\mu}\cdot|a - x|}{\sigma\left(e^{2\mu t} - 1\right)^{3/2}}\cdot\frac 1{\sqrt{2\pi}}\cdot e^{-\frac{\mu(a -
x)^2}{\sigma^2\left(e^{2\mu t} - 1\right)} + 2\mu t} .
\end{equation}

\begin{figure}
\centering
\includegraphics[height=0.33 \textheight]{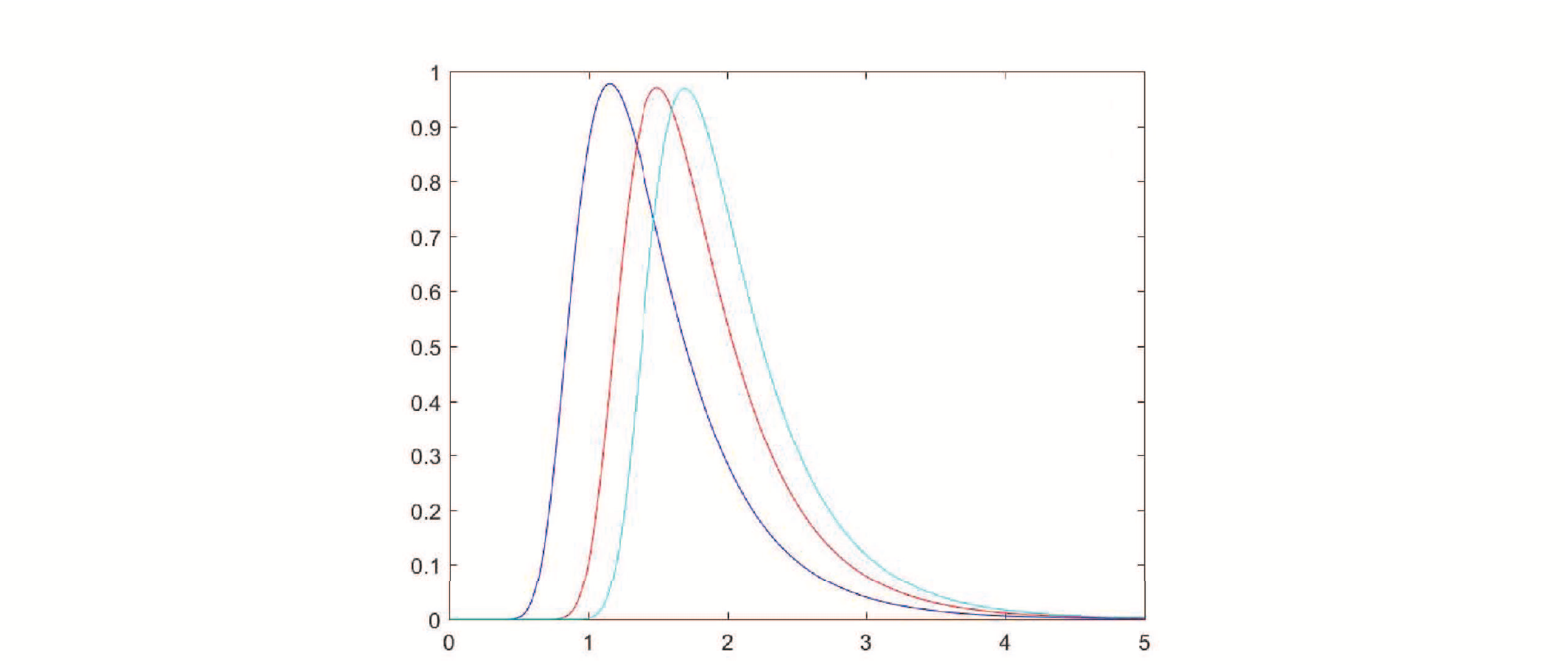}
\caption{Comparison between FPT density of OU through
the boundary $S(t) = ae^{-\mu t}, $ for $x=0, \ \mu = 2, \ \sigma = 0.2$  and
various values of
 $a$; from top to the bottom: $a = 1, a = 2, a = 3 .$  \label{Fig : 9}
}
\end{figure}

In the Figure \ref{Fig : 9} we report the density of $\tau_1(S|x)$, for $x = 0$, $\mu = 2$, $\sigma = 0.2$ and various values
of $a$.

As for the expectation of $\tau_1(a|x),$ we have:
\begin{equation}
E\left[\tau_1(S|x)\right] = \frac{|a - x|\cdot2\mu\sqrt{\mu}}{\sigma\sqrt{\pi}}\cdot\int_0^{+\infty} t  e^{-\frac{(x - a)^2\cdot
\mu}{\sigma^2\cdot\left(e^{2\mu t} - 1\right)} + 2\mu t}\; dt ;
\end{equation}
calculating numerically the integral e.g. for $x=0, \ a=1, \ \mu = \sigma =1,$ one obtains
$E\left[\tau_1(S|x)\right] = 1.025,$
while, Monte Carlo simulation provides the estimate
$\widehat{E}\left[\tau _1(S|x)\right] = 1.048.$ \par
By using \eqref{fT2_TC}, one has:
\begin{align}
f_{T_2}(t) &= \frac{(2\mu)^{5/2}\sigma(a - x)^2}{2\pi^2}\cdot\int_0^{+\infty}\frac{e^{-\frac{\mu(a - x)^2}{\sigma^2(e^{2\mu s} - 1)} + 2\mu(t + 3s)}}{(e^{2\mu(t +
s)} - 1)(e^{2\mu s} - 1)^{3/2}}\cdot\notag \\ &\cdot\left[\int_0^{\frac{\sigma^2}{2\mu}(e^{2\mu(t + s)} - 1)}\frac{e^{-\frac{(a -
x)^2}{2u}}}{u\sqrt{\frac{\sigma^2}{2\mu}(e^{2\mu(t + s)}-1)-u}}du\right]ds . \label{fT2_OU}
\end{align}

\begin{figure}
\centering
\includegraphics[height=0.33 \textheight]{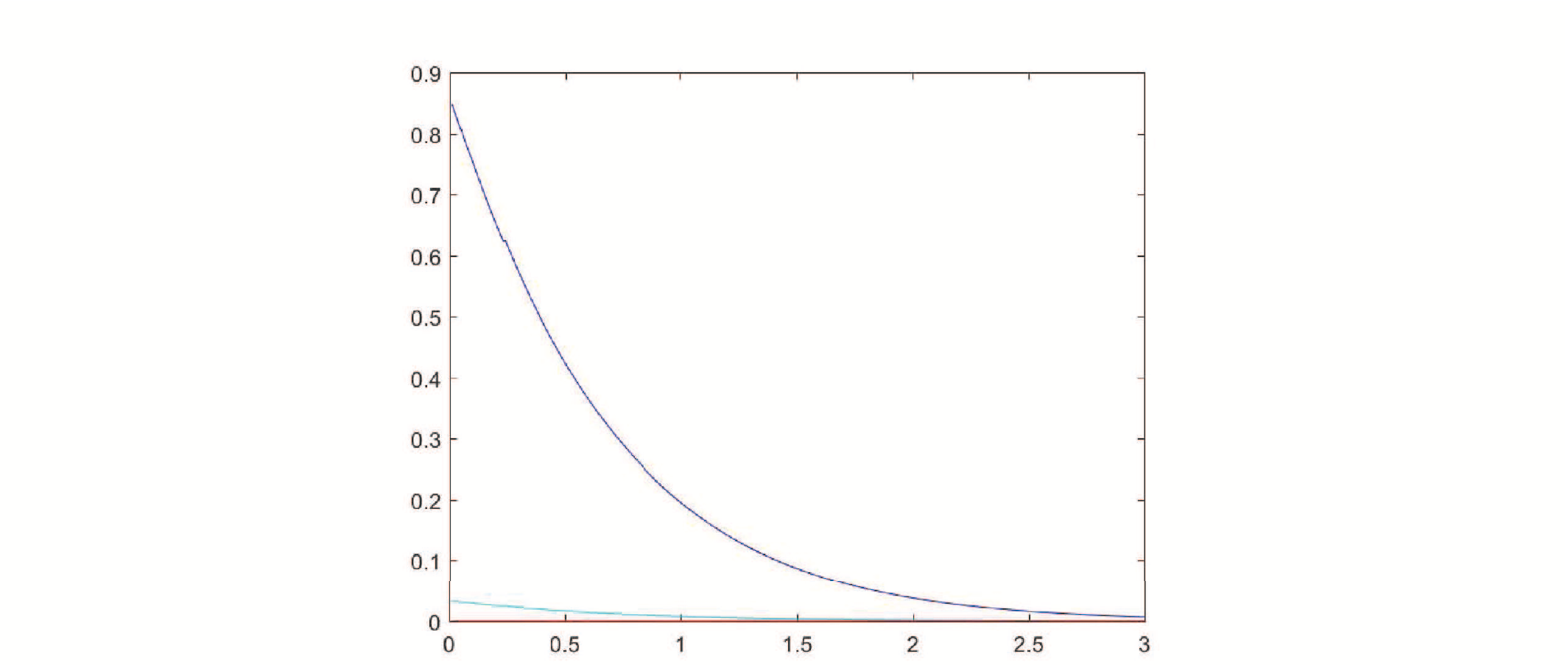}
\caption{Approximate density of $T_2(x)$ for $x = 0, \mu = 2, \sigma = 0.2$ and various values of the parameter
$a;$ from top to the bottom: $a = 1, a = 2, a = 3.$ It refers to the OU process through the boundary $S(t) = ae^{-\mu t}$.
\label{Fig : 10}
}
\end{figure}

\begin{figure}
\centering
\includegraphics[height=0.33 \textheight]{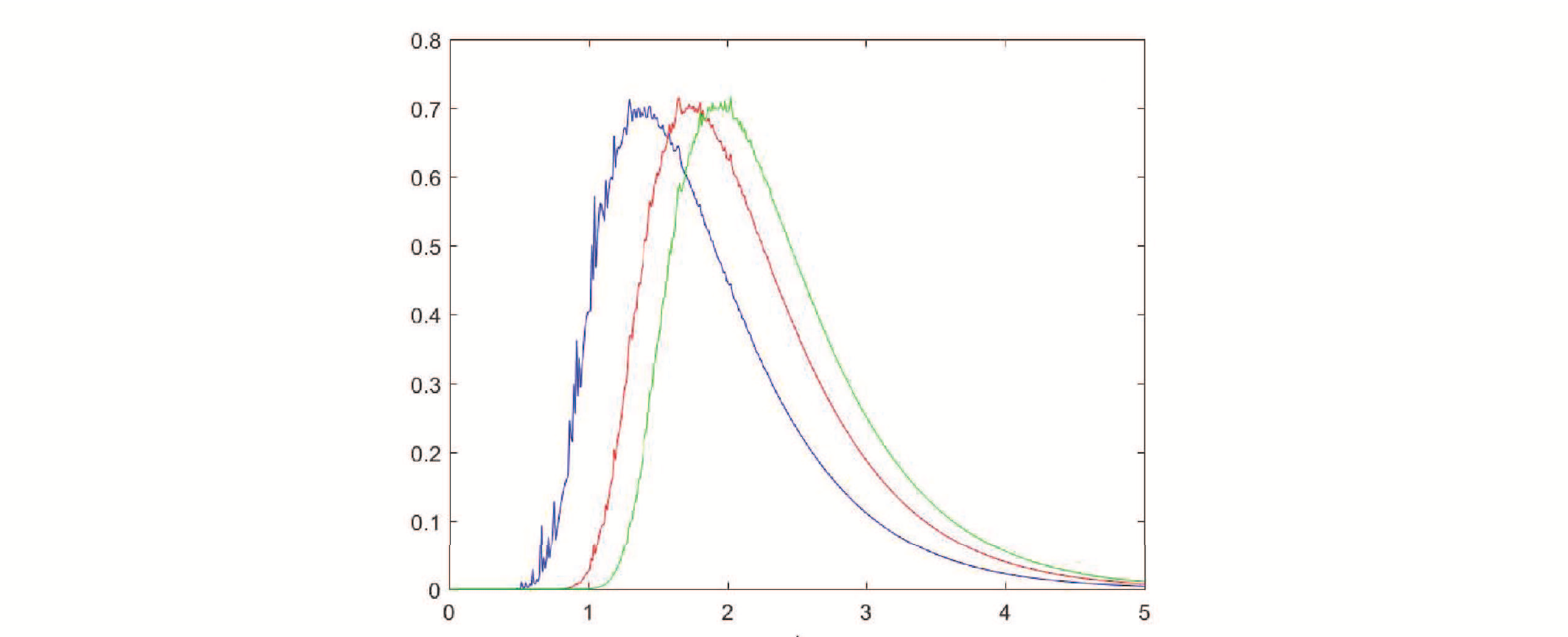}
\caption{Approximate density of $\tau_ 2(S|x)$ for $x = 0, \sigma = 0.2, \mu = 2$ and various values of the
parameter $a$; from top to the bottom: $a = 1, a = 2, a = 3$. It refers to the OU process through the boundary $S(t) = ae^{-\mu t}$.
\label{Fig : 11}
}
\end{figure}

In the Figure \ref{Fig : 10} we report the probability density of $T_2(x)$, obtained from \eqref{fT2_TC}, by calculating numerically the integral, for
$x = 0, \sigma = 0,2, \mu = 2$ and various values of the parameter $a.$ Notice that, for $a=3$ the curve appears to be overlapped
with the time axis.
By using $\eqref{fTau2_TC}$, we obtain
\begin{equation}\label{fTau2_OU}
f_{\tau_2}(t) = \frac{2\mu|a - x|e^{2\mu t}}{\pi\sqrt{2\pi}\cdot(e^{2\mu t} - 1)}\cdot\int_0^{\frac{\sigma^2}{2\mu}(e^{2\mu t} - 1)}\frac{e^{-\frac{(a -
x)^2}{2s}}}{s\sqrt{\frac{\sigma^2}{2\mu}(e^{2\mu t} - 1) - s}}\;ds.
\end{equation}
Calculating numerically the integral, one gets
the density of $\tau_2(S|x);$ for $x = 0, \  \sigma = 0.2,  \ \mu = 2$ and different values of the parameter $a$,
the results are shown in the
Figure \ref{Fig : 11}.

We summarize the LTs of $\tau_1(S|x)$ and $\tau_2(S|x)$ in the case of OU process and $S(t)= ae^{-\mu t}.$ \bigskip

\noindent$ \bullet $ LT of $\tau_1(S|x):$
\begin{equation}
\psi(\lambda) = E\left[e^{- \lambda \tau_1(S|x)}\right] =
\int_0^{+\infty}f_{\tau_1}(t)e^{-\lambda t}\;dt = \frac{(2\mu)^{3/2}\cdot|a - x|}{\sigma\sqrt{2\pi}} \int_0^{+\infty}\frac{e^{2\mu
t - \frac{\mu(a - x)^2}{\sigma^2(e^{2\mu t} - 1)} - \lambda t}}{(e^{2\mu t} - 1)^{3/2}} \ dt .
\end{equation}

$ \bullet $  LT of $\tau_2(S|x):$
$$
\tilde{\psi}(\lambda) = E\left[e^{-\lambda \tau_2(S|x) }\right] = \int_0^{+\infty}f_{\tau_2}(t)e^{-\lambda t}\;dt $$
\begin{equation} \label{LTsecondOU}
=\frac{2\mu|a - x|}{\pi\sqrt{2\pi}}\int_0^{+\infty}\frac{e^{2\mu t - \lambda t}}{e^{2\mu t} - 1}\left[\int_0^{\frac{\sigma^2}{2\mu}(e^{2\mu t} -
1)}\frac{e^{-\frac{(a - x)^2}{2s}}}{s\sqrt{\frac{\sigma^2}{2\mu}(e^{2\mu t} - 1) - s}}\;ds\right]dt .
\end{equation}

\begin{figure}
\centering
\includegraphics[height=0.33 \textheight]{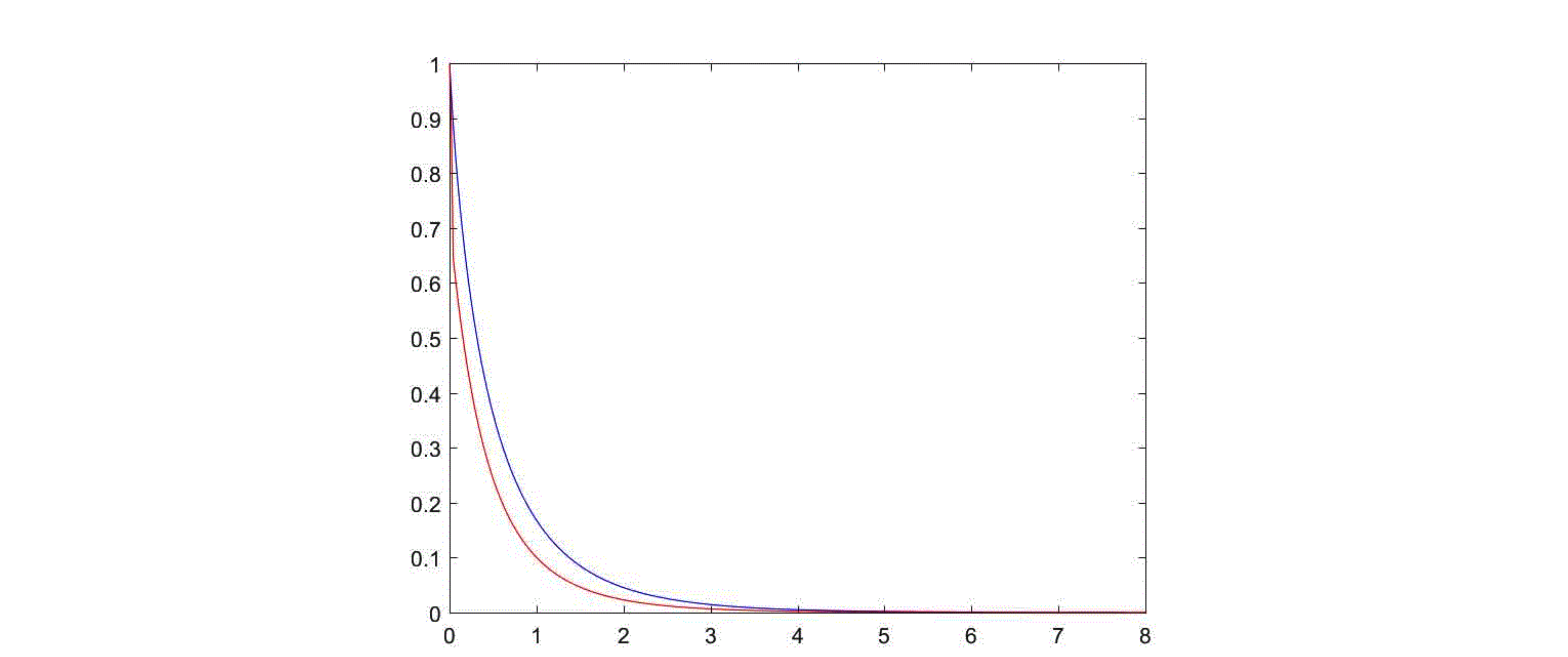}
\caption{Comparison between Laplace transforms of  $\tau_1(S|x)$ (upper curve) and $\tau_2(S|x)$ (lower curve), in the case of OU
process and the boundary $S(t)=ae^ { -\mu t},$ for $\sigma =0.2, \ \mu =2, \ x=0$ and $a=1$. \label{Fig : 15}
}
\end{figure}

As far the expectation of $\tau_2(S|x)$ is concerned, it is
not possible to obtain the exact value,
by taking minus the derivative with respect to $\lambda $  in \eqref{LTsecondOU}, that is:
\begin{equation}
 - \frac {d } {d \lambda } \tilde{\psi}(\lambda)\big{\arrowvert}_{\lambda = 0}= \frac{2\mu|a - x|}{\pi\sqrt{2\pi}}\cdot\int_0^{+\infty}\frac{t
e^{2\mu t}}{e^{2\mu t} - 1}\cdot\left[\int_0^{\frac{\sigma^2}{2\mu}(e^{2\mu t} - 1)}\frac{e^{-\frac{(a - x)^2}{2s}}}{s\sqrt{\frac{\sigma^2}{2\mu}(e^{2\mu t} - 1)
- s}}\;ds\right]\;dt ,
\end{equation}
because of the complexity of the integrals involved in the calculation, so it has to be estimated
by Monte Carlo simulation.
\par
\begin{figure}
\centering
\includegraphics[height=0.33 \textheight]{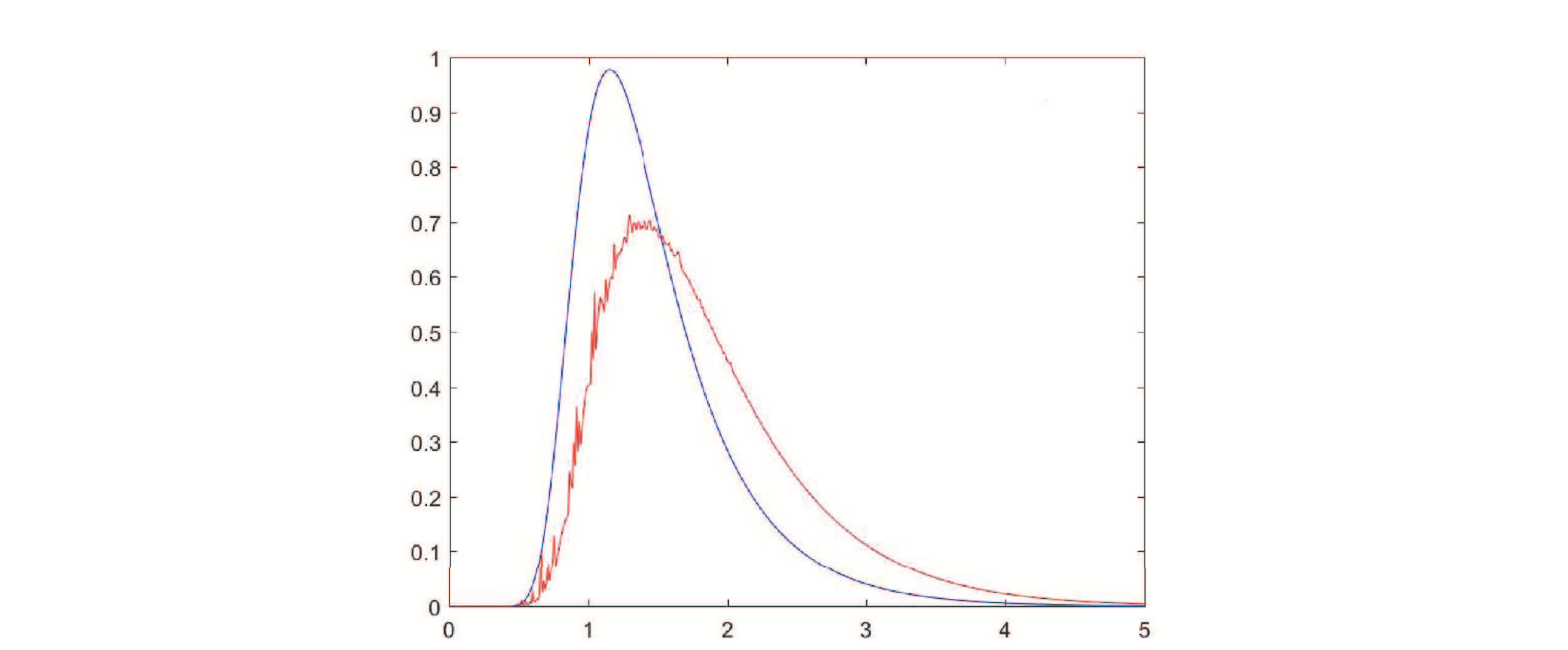}
\caption{Comparison between the density \eqref{fTau2_OU} of $\tau_2(S|x)$ (lower peak) and the
density \eqref{fTau1_OU} of $\tau_1(S|x)$ (upper peak), in the case of OU process and the boundary $S(t)= a e^ { \mu t },$ for $\sigma = 0.2, \mu = 2, \ x=0$ and $a = 1.$ \label{Fig : 12}
}
\end{figure}

In the Figure 14, with regard to OU through the boundary $S(t)= a e^ { - \mu t}, \ a=1,$ we report the comparison between the LT of $\tau_1(S|x)$ and 
that of $\tau_2(S|x),$
for $x = 0, \sigma = 0.2, \mu = 2;$
in the Figure \ref{Fig : 12}, we report the comparison between the density of $\tau_1(S|x)$ and that 
of $\tau_2(S|x)$, for the same set of parameters.

\section{Conclusions and final remarks}
By using the results of \cite{abundo:phy16} on the $n$th-passage time of Brownian motion through  a straight line, we studied the distribution
of the $n$th-passage time through a barrier $a,$ of a diffusion process $X(t),$ given by a space or time transformation of
Brownian motion $B_t;$ precisely, we have considered the cases: $X(t)= v^{-1} (B_t + v(x)),$ where $v(x)$ is an increasing function with $v(0)=0,$
and $X(t)= B ( \rho (t)),$ where $ \rho (t)$  is an increasing function with $\rho (0) =0.$
We also found the Laplace transform of the first and second-passage time of such
processes, and we have reported some explicit examples. \par
Notice that the results of this paper  can be extended to diffusions which are more general than the process $X(t)$  here considered,
for instance to a process of the form
\begin{equation} \label{spacetimetransf}
X(t)= v^ {-1} ( \widehat B(\rho (t)) + v (x) ) ,
\end{equation}
where $\widehat B (t)$ is BM, $v(x)$ and $ \rho (t)$  are regular enough, increasing functions with $\rho (0) = v(0) =0;$ such a process
$X(t)$ is obtained from BM combining a space transformation and a
time-change (see e.g. the discussion in \cite{abundo:stapro12}).
When $\rho (t) = t,$
the process $X(t)$ is conjugated to BM, according to the
definition given in Section 3.1 (see \cite{abundo:stapro12}). The process $X(t)$ given by
\eqref{spacetimetransf} is indeed a weak solution of the SDE:
\begin{equation}
dX(t)= - \frac {\rho ' (t) v'' (X(t)) } { 2 (v' (X(t)))^3 } dt + \frac {\sqrt { \rho '(t)} } { v' (X(t))}  d B_t \ ,
\end{equation}
where  $v'(x)$ and $v''(x)$ denote first and second derivative of
$v(x).$  \par
Provided that the deterministic function $\rho (t)$ is  replaced
with a random function, the representation \eqref{spacetimetransf} is also valid  for a
time homogeneous one-dimensional diffusion driven by the SDE
\begin{equation} \label{eqdiffu}
dX(t) = \mu (X(t)) dt + \sigma (X(t)) d B_t , \ X(0) = x ,
\end{equation}
where the drift $( \mu )$ and diffusion coefficients $ ( \sigma ^2 )$ satisfy the usual conditions
(see e.g. \cite{ikeda:81})
for existence and uniqueness
of the solution of \eqref{eqdiffu}. In fact,
let $w(x) $ be the {\it scale function} associated to
the diffusion $X(t)$ driven by the SDE \eqref{eqdiffu}, that is,
the  solution of $L w(x) =0, \ w(0) =0, \ w'(0) =1 ,$
where $L$ is the infinitesimal generator of $X$ given by
$
L h = \frac 12 \sigma ^2 (x) \frac {d^2 h } {d x^2}  + \mu (x) \frac {d h } {d x}.
$
As easily seen,  if the integral $ \int _0 ^t \frac { 2 \mu (z)}
{\sigma ^2 (z)} \ dz $ converges, the scale function is explicitly given by
\begin{equation} \label{scalefunction}
w(x)= \int _ 0 ^x \exp \left ( - \int _0  ^t \frac { 2 \mu (z)}  {\sigma ^2 (z)} \ dz \right ) dt .
\end{equation}
If $\zeta (t) := w(X(t)),$ by It${\rm \hat o}$'s formula one
obtains
\begin{equation}
 \zeta (t) = w ( x) + \int _0 ^t w '( w ^{-1} (\zeta (s)) ) \sigma ( w ^{-1} (\zeta  (s))) d B_s  \ ,
\end{equation}
that is, the process $\zeta (t) $ is a local martingale, whose
quadratic variation is
\begin{equation} \label{rho}
\rho (t) \doteq \langle \zeta  \rangle _t = \int _0 ^ t [ w' (X(s)) \sigma (X(s)) ] ^ 2 ds, \ t \ge 0 .
\end{equation}
The (random) function $\rho(t)$ is differentiable, increasing, and
$\rho (0)=0.$ If $\rho ( + \infty) = + \infty,$
by the Dambis, Dubins-Schwarz theorem (see e.g. \cite{revuzyor:con91}) one gets
that there exists a BM $\widehat
B$ such that $\zeta (t) = \widehat B ( \rho (t)) + w ( \eta ).$
Thus, since $w$ is invertible, one obtains the representation \eqref{spacetimetransf} with $w$ in place of $u.$ \par\noindent
Notice, however, that the successive-passage time problem for the process $X(t)$ given by \eqref{spacetimetransf} can be
addressed as in the  case when $\rho (t) $ is a deterministic function, only if
the distribution of the random process $\rho (t)$ is explicitly known.


\begin{thebibliography}{99}

\bibitem {abundo:phy16}
Abundo, M., 2016. \newblock
On the excursions of drifted Brownian motion and the successive passage times of Brownian motion.
\newblock{Physica A} 457, 176–182. doi:10.1016/j.physa.2016.03.052

\bibitem {abundo:smj13}
Abundo, M., 2013. \newblock
On the representation of an integrated Gauss-Markov process.
\newblock{Sci. Math. Jpn. Online} e-2013, 719–-723.

\bibitem {abundo:smj15}
Abundo, M., 2015. \newblock
On the first-passage time of an integrated Gauss-Markov process.
\newblock{Sci. Math. Jpn. Online} e-2015, 28,  1--14.



\bibitem {abundo:stapro12}
Abundo, M. 2012.\newblock
An inverse first-passage problem for one-dimensional diffusions with random
      starting point.
\newblock{\it Statist. Probab. Lett.} 82(1):7–-14.

\bibitem  {abundo:ric01}
Abundo, M. 2001.\newblock
Some results about boundary crossing for Brownian motion.
\newblock{\it Ricerche di Matematica} vol. L(2): 283--301.



\bibitem  {ell:2000}
Elliot, R.J., Jeanblanc, M., and Yor, M., 2000.
\newblock
On models of default risk. \newblock
{Math. Finance}, 10, 179--196.

\bibitem{ikeda:81}
Ikeda, N. and Watanabe, S., 1981. \newblock Stochastic differential
equations and diffusion processes. North-Holland Publishing
Company.

\bibitem  {jea:2009}
Jeanblanc, M., Yor, M., and Chesney, M., 2009. \newblock
Mathematical Methods for financial Markets. \newblock
Springer Finance, Springer.

\bibitem  {kar:98}
Karatzas, I., Shreve, S., 1998. \newblock
Brownian Motion and Stochastic Calculus. \newblock
Springer.

\bibitem  {kleb:05}
Klebaner, F.C., 2005. \newblock
Introduction to Stochastic Calculus with Applications (Second Edition). \newblock
Imperial College Press, London.


\bibitem
{revuzyor:con91}
Revuz, D. and Yor, M., 1991. \newblock
Continous martingales and Brownian motion. \newblock
Springer-Verlag, Berlin Heidelberg.

\bibitem
{ross:pmodel10}
Ross, S.M., 2010. \newblock
Introduction to Probability Models. Tenth Edition. \newblock
Academic Press, Elsevier, Burlington.

\bibitem
{salm:88}
Salminen, P., 1988. \newblock
On the First Hitting Time and the Last Exit Time for a Brownian Motion to/from a Moving Boundary. \newblock
Adv. Appl. Probab.
20 (2), pp. 411--426.

\bibitem
{sch:39}
Schumpeter, J.A., 1939. \newblock
Business Cycles. A Theoretical, Historical and Statistical Analysis of the Capitalist
Process. \newblock
Mc Graw-Hill, New York.


\end{thebibliography}
\end{document}